\documentclass[a4paper,12pt,reqno]{amsart}
\usepackage{amsmath, amsthm, amssymb, color}
\usepackage[colorlinks=true,linkcolor=blue,urlcolor=blue,citecolor=teal]{hyperref}
\usepackage{graphicx}
\usepackage{caption,soul}
\usepackage{mathtools}
\usepackage{verbatim}
\usepackage{tikz}
\usepackage[all]{xy}
\usepackage[english]{babel}
\usepackage[letterpaper,top=2cm,bottom=2cm,left=2cm,right=2cm]{geometry}
\usepackage{cleveref} 
\usepackage{graphicx}
\usepackage{nicematrix}
\usepackage{wasysym}
\usepackage{tikz-cd}
\usepackage{enumitem}
\usetikzlibrary{arrows}
\usepackage{algorithm}
\usepackage{shuffle}
\usepackage[normalem]{ulem}
\usepackage{subcaption}

\usepackage{courier} 

\newtheorem{theorem}{Theorem}[section]

\newtheorem{lemma}[theorem]{Lemma}
\newtheorem{corollary}[theorem]{Corollary}
\theoremstyle{definition}
\newtheorem{definition}[theorem]{Definition}
\newtheorem{remark}[theorem]{Remark}

\newtheorem{example}[theorem]{Example}

\newcommand{\RR}{\mathbb{R}}

\newcommand{\CC}{\mathcal{C}}
\newcommand{\E}{\mathcal{E}}

\newcommand{\TT}{\mathbb{T}}

\newcommand{\shu}{\shuffle}

\newcommand{\x}{\mathbf{x}}
\newcommand{\p}{\mathbf{p}}
\newcommand{\q}{\mathbf{q}}
\newcommand{\y}{\mathbf{y}}

\newcommand{\vv}{\mathbf{v}}

\newcommand{\N}{\mathbb{N}}
\newcommand{\R}{\mathbb{R}}

\newcommand{\ot}{\otimes}

\newcommand{\id}{\mathrm{id}}

\newcommand{\SL}{\mathop{\rm SL}\nolimits}

\newcommand{\rk}{\mathop{\rm rk}\nolimits}

\newcommand{\im}{\mathop{\rm im}\nolimits}

\newcommand{\f}{\varphi}

\def\word#1{{\texttt{#1}}}

\newcommand{\PW}[1]{\mathcal{C}^{{#1}\,\textrm{-pw}}}

\title{
Path signatures of ODE solutions}

\author[F. Galuppi]{Francesco Galuppi}
\address[F. Galuppi]{Faculty of Mathematics, Informatics, and Mechanics, University of Warsaw, Banacha 2, 02-097 Warsaw, Poland}
\thanks{Galuppi acknowledges support from the National Science Center, Poland, under the project ``Tensor rank and its applications to signature tensors of paths'', 2023/51/D/ST1/02363.}
\email{galuppi@mimuw.edu.pl (ORCID 0000-0001-5630-5389)}

\author[G. Moreno]{Giovanni Moreno}
\address[G. Moreno]{Faculty of Physics, University of Warsaw, Pasteura 5, 02-093
Warsaw, Poland}
\email{gmoreno@fuw.edu.pl (ORCID 0000-0001-8670-724X)}

\author[P. Santarsiero]{Pierpaola Santarsiero}
\address[P. Santarsiero]{Universit\`a di Bologna, Dipartimento di Matematica, Piazza di Porta S. Donato 5, 40126 Bologna, Italy}
\email{pierpaola.santarsiero@unibo.it (ORCID 0000-0003-1322-8752)}
\thanks{Santarsiero was supported by the European Union under NextGenerationEU. PRIN 2022, Prot. 2022E2Z4AK}

\begin{document}

\begin{abstract}
The signature of a path is a sequence of tensors which allows to uniquely reconstruct the path. By employing the geometric theory of nonlinear systems of ordinary differential equations, we find necessary and sufficient algebraic conditions on the signature tensors of a path to be a solution of a given system of ODEs. As an application, we describe in detail the systems of ODEs that describe the trajectories of a vector field, in particular a linear and Hamiltonian one. 
\end{abstract}

\maketitle

\section{Introduction}\label{section: intro}

The signature of a path $X:[a,b]\to\R^d$ is a sequence of tensors, introduced in \cite{Chenoriginal}, that encodes the essential information contained in $X$. Under the assumption that $X$ is sufficiently smooth, the $k$-th signature of $X$ is the tensor $\sigma^{(k)}(X)\in (\R^d)^{\ot k}$ whose $(i_1,\dots,i_k)$-th entry is the iterated integral
\begin{equation}\label{defin di signature con gli iterated integrals}
\sigma^{(k)}(X)_{i_1,\dots,i_k}=\int_{a}^{b}\int_{a}^{t_k}\dots\int_{a}^{t_3}\int_{a}^{t_2}\dot{X}_{i_1}(t_1)\cdots\dot{X}_{i_k}(t_k){\rm d}t_1\cdots {\rm d}t_k.
\end{equation}
By convention we set $\sigma^{(0)}(X)=1$. The \emph{signature} of $X$ is the sequence $\sigma(X)=(\sigma^{(k)}(X)\mid k\in\N)$. As proven in \cite{chenuniqueness}, the signature $\sigma(X)$ allows to uniquely reconstruct the path $X$, up to a mild equivalence relation. 
After Chen's pioneering work, signatures proved to be a powerful tool in many areas of mathematics, notably in stochastic analysis, and in particular in the theory of rough paths \cite{lyons2014rough}. The twenty-first century saw an explosion of applications of signature tensors: from topological data analysis \cite{signaturesetopdata}, where signatures are used to extract features of barcodes, to quantum field theory \cite{quantumfield}, where they are called Dyson series and have deep relationships with conformal field theory and Feynman diagram computations. Signatures are successfully employed also in medical statistics, for instance to facilitate the formulation of diagnosis \cite{medical}. The theory keeps developing and finding new applications, for instance in machine learning \cite{machine} and cybersecurity \cite{cyber}, just to mention a few.

In this paper we approach signature tensors from a geometric viewpoint.
Given a hyperplane $H$, \cite{GalSan} establishes  a necessary and sufficient condition on $
\sigma(X)$ in order for $X$ to lie on $H$. It is then natural  to study paths contained in an arbitrary algebraic variety: we develop this idea and we apply it to the characterization of solutions of systems of ordinary differential equations (ODE) in terms of signature tensors.\par

In \Cref{section: paths su una varietà} we show how to use signature tensors to characterize paths lying on a given algebraic variety. \Cref{section: lagrangian paths} takes inspiration from differential geometry: we give necessary and sufficient algebraic conditions on $\sigma(X)$ in order for $X$ to be holonomic. Combining these two results, in \Cref{section: ODE} we characterize solutions of Cauchy problems in terms of the signature of a suitable path. In \Cref{section: esempi e applicazioni} the main result is  applied to some family of paths that are interesting in differential geometry and mathematical physics, such as integral and Hamiltonian paths.\par 

\subsubsection*{Acknowledgements} Galuppi and Moreno are grateful to Jarosław Buczy\'nski for the seminar ``Progedia'', where the idea of this paper was born. The authors thank Rosa Preiß for interesting and helpful discussions. Moreover we tank Terry Lyons and Peter Friz for their stimulating questions.

\section{Chen's uniqueness theorem and the shuffle identity}
In this section we fix the notations we will use through the paper, starting from the class of paths we will be working with. In \cite{chenuniqueness}, Chen deals with piecewise regular paths, which form a subclass of $\PW{1}([a,b])$. Definition \eqref{defin di signature con gli iterated integrals} applies to a larger class, namely bounded variation paths. The theory has been widely extended and generalized, so now it is even possible to consider the signature of a rough path, as illustrated in \cite{lyonsdriven}. Since we are interested in characterizing solutions of a Cauchy problem, in this paper we assume that the path $X$ is continuous and piecewise differentiable.

\begin{definition}\label{definition: paths of class Ck}
Let $[a,b]\subseteq\R$ be an interval and let $k\in\N$. A continuous path $X:[a,b]\to\R^d$ is \emph{piecewisely of class $\CC^k$} if 
there exist $a=t_0<t_1<\dots<t_{n-1}<t_n=b$ such that the restriction $X|_{[t_{i-1},t_i]}$ is of class $\CC^k([t_{i-1},t_i])$ for every $i\in\{1,\ldots,n\}$. The class of continuous paths piecewisely of class $\CC^k$ is denoted by $\PW{k}([a,b])$.
\end{definition}
To uniquely recover the path from its signature, often we will also assume that $X$ has no tree-like excursions.
\begin{definition}
Let $X:[a,b]\to\R^d$ be a path. The \emph{inverse path} of $X$, denoted by $X^{-1}$, is the path $X^{-1}:[a,b]\to\R^d$ obtained by changing the orientation of $X$, namely $$X^{-1}(t)=X\left(\frac{(t-a)a}{b-a}+\frac{(b-t)b}{b-a}\right).$$
The path obtained by concatenating $X$ and $X^{-1}$ is called a \emph{tree-like excursion}. If a path $X$ contains no tree-like excursion, then it is called \emph{reduced}. Notice that in \cite[Definition 3.1]{chenuniqueness} the same property is called being \emph{irreducible}. Denote by $X_{red}$ the unique reduced path which coincides with $X$ up to tree-like excursions.
\end{definition}

With the notion of reduced path, we can precisely state the uniqueness theorem.

\begin{theorem}\label{theorem: chen uniqueness}
Let $X,Y\in\PW{1}([a,b])$ be paths. Then $\sigma(X)=\sigma(Y)$ if and only if $X_{red}$ and $Y_{red}$ are the same path, up to translation.
\end{theorem}

\Cref{theorem: chen uniqueness} appeared in \cite[Theorem 4.1]{chenuniqueness} and has been generalized in \cite[Theorem 4]{HL10} and \cite[Theorem 1.1]{uniquenessrough} to more general classes of paths. We conclude this section by recalling the interpretation of signature tensors in terms of combinatorics and non-commutative algebra. The textbook reference for this material is \cite{Reutenauer}.
\begin{definition}\label{definition: word and shuffle}
Let $d\in\N$. We regard every element of $\{1,2,\dots,d\}$ as a \emph{letter}. We denote letters by $\word{1},\word{2},\dots, \word{d}$ in order to highlight the difference between the number 1 and the letter $\word{1}$. The set $\{\word{1},\word{2},\dots, \word{d}\}$ is called  \emph{alphabet}. An ordered sequence of letters $\word{i}_1 \word{i}_2 \cdots \word{i}_k$ is called a \emph{word} of length $k$. The only word of length zero is called the \emph{empty word} and denoted by $\varnothing$. The \emph{concatenation} of two words $w_1$ and $w_2$ is the word $w_1 w_2$ obtained by writing $w_1$ followed by $w_2$, and it is extended by linearity.
\end{definition}

We denote by $\TT(\R^d)$ the graded $\R$-vector space whose basis consists of all words in the alphabet $\{\word{1},\dots,\word{d}\}$. Elements of $\TT(\R^d)$ can be considered linear functions on the set of signatures. Indeed, if $X:[a,b]\to\R^d$ is a path, then
\begin{equation*}
\langle\sigma(X),\word{i}_1 \word{i}_2 \cdots \word{i}_k\rangle=\sigma^{(k)}(X)_{i_1,i_2,\dots,i_k}
\end{equation*}
is defined to be the $(i_1,i_2,\dots,i_k)$-th entry of $\sigma^{(k)}(X)$. This definition is extended by linearity, by setting
\begin{equation*}
\langle\sigma(X),w_1+w_2\rangle=\langle\sigma(X),w_1\rangle+\langle\sigma(X),w_2\rangle.
\end{equation*}
The set of all words   is endowed with  a non-associative, non-commutative operation that will be useful for our arguments in \Cref{section: paths su una varietà}.

\begin{definition}\label{definition: half shuffle}
The \emph{right half-shuffle} between two words, denoted by $\succ$, is defined recursively by
\begin{align*}
w_1\succ \word{i}&=w_1\word{i}\mbox{ for every letter } \word{i}\in \{\word{1},\word{2},\dots, \word{d}\}\\
w_1\succ w_2\word{i}&=(w_1\succ w_2+w_2\succ w_1)\word{i}.
\end{align*}
The \emph{shuffle} between $w_1$ and $w_2$ is $w_1\shu w_2=w_1\succ w_2+w_2\succ w_1$.
\end{definition}

For example, $\word{1}\succ\word{23}=\word{123}+\word{213}$ and $\word{1}\shu\word{23}=\word{123}+\word{213}+\word{231}$. A reference for the half-shuffle is \cite{halfshuffle}. Unlike the half-shuffle, the shuffle is associative and commutative, so $(\TT(\R^d),\shu,\varnothing)$ is a graded commutative $\R$-algebra. The connection between shuffle product and signature tensors is the shuffle identity, that appears for instance in \cite[Proof of Corollary~3.5]{Reutenauer}.
\begin{lemma}\label{lem:shuffle identity}
Let $X:[a,b]\to \R^d$ is a path. If $w_1$ and $w_2$ are words in the alphabet $\{\word{1},\dots,\word{d}\}$, then $\langle\sigma(X),w_1\shu w_2\rangle=\langle\sigma(X),w_1\rangle\cdot \langle\sigma(X),w_2\rangle$.
\end{lemma}

\section{Signatures of paths on an algebraic variety}\label{section: paths su una varietà}
Given a polynomial $g\in\R[x_1,\dots,x_d]$, the goal of this section is to characterize the paths in $\R^d$ whose images lie on the affine hypersurface defined by $g$. Our strategy is to apply a suitable polynomial map $p:\R^d\to\R^{d+1}$ and track down how signatures are transformed under $p$, according to \cite{CP2020}.  
\begin{definition}\label{definition: p and M and phi and g}
Let $g\in\R[x_1,\dots,x_d]$ and consider the map $p:\R^d\to\R^{d+1}$ defined by
\[p(x_1,\dots,x_d)=(x_1,\dots,x_d,g(x_1,\dots,x_d)).
\]
We will also  need the ring homomorphism $\f_d:\R[x_1,\dots,x_d]\to (\TT(\R^d),\shu,\varnothing)$ given by $\f_d(x_i)= \word{i}$.
\end{definition}

Throughout this section   $g\in\R[x_1,\dots,x_d]$ and $p:\R^d\to\R^{d+1}$  will be as  in \Cref{definition: p and M and phi and g}. The following  statement appears  in \cite[Theorem 1 and Corollary 1]{CP2020}.

\begin{theorem}\label{theorem: Laura e Rosa} Let 
$X:[a,b]\to\R^d$ be a path, define $\widetilde{p}(x)=p(x+X(a))-p(X(a))$, and let $J$ be the $(d+1)\times d$ Jacobian matrix of $\widetilde{p}$. Then the map $M_{\widetilde{p}}: (\TT(\R^{d+1}),\shu,\varnothing)\to (\TT(\R^d),\shu,\varnothing)$, recursively defined on words by
\begin{align*}
M_{\widetilde{p}}&(\varnothing)=\varnothing ,\\
M_{\widetilde{p}}&(v  \word{i})=\sum_{j=1}^d (M_{\widetilde{p}}(v)\shu \f_d(J_{i,j}))\word{j} ,
\end{align*}
is both an $\R$-algebra homomorphism and a half-shuffle homomorphism. Moreover 
\[
\langle\sigma(p\circ X),v\rangle=\langle\sigma(X),M_{\widetilde{p}}(v)\rangle \mbox{ for every } v\in\TT(\R^{d+1}).\qedhere
\]
\end{theorem}
Before we move to the main theorem, we state and prove some small technical statement.  The next two results will help us understand how $M_{\widetilde{p}}$ works.
\begin{lemma}\label{lemma: proprietà M}
Let 
$X:[a,b]\to\R^d$ be a path and let $\widetilde{p}(x)=p(x+X(a))-p(X(a))$.
\begin{enumerate}
\item \label{item: M identity on letters} If $\word{i}\in\{\word{1},\dots,\word{d}\}$, then $M_{\widetilde{p}}(\word{i})=\word{i}$.
\item \label{item: M assorbe lettere} If $v\in\TT(\R^{d+1})$ and $\word{i}\in\{\word{1},\dots,\word{d}\}$, then $M_{\widetilde{p}}(v\word{i})=M_{\widetilde{p}}(v)\word{i}$.
\item \label{item: le nostre equazioni sono M(ep)}If $w\in\TT(\R^d)$, then $$M_{\widetilde{p}}((\word{d+1})w)=\sum_{j=1}^d\f_d\left(\frac{\partial}{\partial x_j}g(x+X(a))\right)\word{j} w.$$
\end{enumerate}
\begin{proof}
First observe that the Jacobian matrix of $\widetilde{p}$ is
\[J=\begin{pmatrix}
&&\\&&\\
&\id_{d\times d} &\\
&&\\&&\\
\hline\\[-2ex]
\frac{\partial}{\partial x_1}g(x+X(a)) &\dots &\frac{\partial}{\partial x_d}g(x+X(a))
\end{pmatrix}.\]
If $\word{i}\in\{\word{1},\dots,\word{d}\}$, then by definition
\begin{align*}
   M_{\widetilde{p}}(\word{i})
=  \sum_{j=1}^d( M_{\widetilde{p}}(\varnothing)\shuffle\f_d( J_{i,j}))\word{j}
=  \sum_{j=1}^d( \varnothing\shuffle\delta_{ij}\varnothing)\word{j}=  
\word{i},
\end{align*}
so claim  \eqref{item: M identity on letters} holds. Now, if $v\in\TT(\R^{d+1})$ and $\word{i}\in\{\word{1},\dots,\word{d}\}$, then the fact that $M_{\widetilde{p}}$ is a $\succ$-morphism implies that
\[M_{\widetilde{p}}(v\word{i})=M_{\widetilde{p}}(v\succ \word{i})=M_{\widetilde{p}}(v)\succ M_{\widetilde{p}}(\word{i}) =M_{\widetilde{p}}(v)\word{i},\]
where the last equality follows from \eqref{item: M identity on letters}, so claim  \eqref{item: M assorbe lettere} holds. To  prove claim  \eqref{item: le nostre equazioni sono M(ep)} assume first that   $w$ is the empty word. Then 
\begin{align*}
M_{\widetilde{p}}((\word{d+1})\varnothing)=M_{\widetilde{p}}(\word{d+1})= \sum_{j=1}^d( M_{\widetilde{p}}(\varnothing)\shuffle \f_d(J_{d+1, j}))\word{j}= \sum_{j=1}^d \f_d\left(\frac{\partial}{\partial x_j}g(x+X(a))\right)\word{j}.
\end{align*}
In order to conclude, let  $w$ be  any word in the alphabet $\{\word{1},\dots,\word{d}\}$. It suffices then to   concatenate the letters of $w$ one by one on the right and apply  \eqref{item: M assorbe lettere}.
\end{proof}\end{lemma}

The definition of $M_{\widetilde{p}}$ and   \Cref{lemma: proprietà M} show that the $(d+1)$-st coordinate behaves in a different way than the others. Next lemma allows us to reduce the case of words containing multiple instances of the letter $\word{d+1}$ to the case of  words containing  only one instance.

\begin{lemma}\label{lemma: da avere tante e ad avere una sola e}
Let 
$X:[a,b]\to\R^d$ be a path and define  $\widetilde{p}(x)=p(x+X(a))-p(X(a))$. If $v\in\TT(\R^{d+1})$ is a word containing the letter $\word{d+1}$, then there exist words $v_1,\dots,v_t\in\TT(\R^{d+1})$, each containing $\word{d+1}$ only once, such that $M_{\widetilde{p}}(v)$ is a linear combination of $M_{\widetilde{p}}(v_1),\dots,M_{\widetilde{p}}(v_t)$.
\begin{proof}
Let $n$ be the number of times the letter $\word{d+1}$ appears in $v$. We argue by induction on $n$. If $n=1$, then it suffices to take $t=1$ and $v_1=v$. Now assume that $n\ge 2$ and write $v$ as the concatenation $v=u\cdot(\word{d+1})\cdot q$, where $q$ is a word not containing $\word{d+1}$. Then
\begin{align*}
M_{\widetilde{p}}(u\cdot(\word{d+1}))&=\sum_{j=1}^d(M_{\widetilde{p}}(u)\shu \f_d(J_{d+1,j}))\word{j}
= \sum_{j=1}^d(M_{\widetilde{p}}(u)\shu 2(\word{j}+X_j(a)\varnothing))\word{j}\\&= 2\sum_{j=1}^d(M_{\widetilde{p}}(u)\shu \word{j})\word{j}+2\sum_{j=1}^d X_j(a)(M_{\widetilde{p}}(u)\shu\varnothing))\word{j}
\\&= 2\sum_{j=1}^d(M_{\widetilde{p}}(u)\shu M_{\widetilde{p}}(\word{j}))\succ M_{\widetilde{p}}(\word{j})+2\sum_{j=1}^d X_j(a)M_{\widetilde{p}}(u)\succ M_{\widetilde{p}}(\word{j})\mbox{ by \Cref{lemma: proprietà M}}
\\&= 2\sum_{j=1}^d(M_{\widetilde{p}}(u\shu\word{j}))\succ M_{\widetilde{p}}(\word{j})+2\sum_{j=1}^d X_j(a)M_{\widetilde{p}}(u\word{j})\mbox{ by \Cref{theorem: Laura e Rosa}.}
\end{align*}
By construction $u$ contains the letter $\word{d+1}$ only $n-1$ times, so $u\word{j}$ and every summand of $u\shu\word{j}$ also contains $\word{d+1}$ only $n-1$ times. By induction hypothesis, we can decompose $M_{\widetilde{p}}(u\word{j})$ and every summand of $M_{\widetilde{p}}(u\shu\word{j})$. In other words there exist $u_1,\dots,u_s\in\TT(\R^{d+1})$, each containing the letter $\word{d+1}$ only once, such that
\begin{align*}
M_{\widetilde{p}}(u\cdot(\word{d+1}))=\sum_{j=1}^d\sum_{i=1}^s M_{\widetilde{p}}(u_i)\succ M_{\widetilde{p}}(\word{j})=\sum_{j=1}^d\sum_{i=1}^s M_{\widetilde{p}}(u_i\word{j}).
\end{align*}
To finish the proof, we write $q=\word{i}_1\word{i}_2\dots \word{i}_k$ as a concatenation of letters in the alphabet $\{\word{1},\dots,\word{d}\}$. Since $M_{\widetilde{p}}$ is a $\succ$-homomorhpism,  \Cref{lemma: proprietà M} implies that  
\begin{align*}M_{\widetilde{p}}(u\cdot(\word{d+1})\cdot \word{i}_1)&=\sum_{j=1}^d\sum_{i=1}^s M_{\widetilde{p}}(u_i\word{j})\succ \word{i}_1=\sum_{j=1}^d\sum_{i=1}^s M_{\widetilde{p}}(u_i\word{j})\succ M_{\widetilde{p}}(\word{i}_1)
\\&=\sum_{j=1}^d\sum_{i=1}^s M_{\widetilde{p}}(u_i\word{j}\succ \word{i}_1)=\sum_{j=1}^d\sum_{i=1}^s M_{\widetilde{p}}(u_i\word{j} \word{i}_1).\end{align*}
By concatenating one by one all the letters of $q$, we obtain 
\[M_{\widetilde{p}}(v)=M_{\widetilde{p}}(u\cdot(\word{d+1})\cdot \word{i}_1\word{i}_2\dots\word{i}_k)=\sum_{j=1}^d\sum_{i=1}^s M_{\widetilde{p}}(u_i\word{j} \word{i}_1\word{i}_2\dots\word{i}_k).\qedhere
\]
\end{proof}
\end{lemma}

One last simplification is in order:  reduce the case of words \textit{containing} the letter $\word{d+1}$ to the case of words \textit{beginning} with the letter $\word{d+1}$.

\begin{lemma}\label{lemma: ridurre da una d+1 all'inizio ad una d+1 ovunque}
Let $\word{a} $ be a letter and let $w$ be a word containing $\word{a}$. Then there exists words $q_1,\dots,q_t$ and $r_1,\dots,r_t$ such that
\begin{equation}\label{equation: appartiene a Ishu(M(W))}w\mbox{ is a linear combination of }
\word{a} q_1\shuffle r_1,\dots,\word{a} q_t\shuffle r_t.
\end{equation}
\begin{proof}
Let $n$ be the position of the first appearence of $\word{a}$ in $w$. We argue by induction on $n$. If $n=1$, then we conclude by writing $w=\word{a} q_1$, taking $t=1$ and $r_1=\varnothing$.
If $n\ge 2$, then we can write $w=w_0\word{a} w_1$, where $w_0$ does not contain the letter $\word{a}$. By \Cref{lem:shuffle identity}, 
\[w_0\shu \word{a} w_1=w_0\word{a} w_1+\sum_{h=1}^m u_h,
\]
where each $u_h$ is a word where $\word{a}$ appears in position $n-1$ or earlier. By induction hypothesis, each $u_h$ admits a decomposition  as in  \eqref{equation: appartiene a Ishu(M(W))}, hence
\[w=w_0\word{a} w_1=\word{a} w_1\shu w_0-\sum_{h=1}^m u_h
\]
admits such a decomposition as well. 
\end{proof}
\end{lemma}

Combining \Cref{lemma: da avere tante e ad avere una sola e} and \Cref{lemma: ridurre da una d+1 all'inizio ad una d+1 ovunque}, we obtain the technical tool that we will need.

\begin{corollary}\label{corollary: any words is a shuffle combination of stuff with only one e at the beginning}
Let 
$X:[a,b]\to\R^d$ be a path and define  $\widetilde{p}(x)=p(x+X(a))-p(X(a))$. If $v\in\TT(\R^{d+1})$ is a word containing the letter $\word{d+1}$, then there exist $q_1,\dots,q_h,r_1\dots,r_h\in\TT(\R^{d+1})$ such that the letter $\word{d+1}$ does not appear in any of the $q_i$'s and $$M_{\widetilde{p}}(v)
=\sum_{i=1}^h M_{\widetilde{p}}((\word{d+1})q_i)\shu M_{\widetilde{p}}(r_i)
.$$
\end{corollary}

\begin{example}We illustrate two instances of   \Cref{corollary: any words is a shuffle combination of stuff with only one e at the beginning}. For $d=2$ and $v=\word{1}\word{3}\in\TT(\R^{d+1})$, we have
$M_{\widetilde{p}}(\word{1}\word{3})=M_{\widetilde{p}}(\word{3})\shu M_{\widetilde{p}}(\word{1})-M_{\widetilde{p}}(\word{3}\word{1})$. Instead, if we take $v=\word{3}\word{3}\in\TT(\R^{d+1})$, then
\begin{align*}
M_{\widetilde{p}}(\word{3}\word{3})&=2\sum_{j=1}^2(M_{\widetilde{p}}(\word{3})\shu (\word{j}+X_j(a)\varnothing))\word{j}
=2\sum_{j=1}^2(M_{\widetilde{p}}(\word{3})\shu \word{j})\word{j}+2\sum_{j=1}^2X_j(a)M_{\widetilde{p}}(\word{3})\word{j}
\\&=2\sum_{j=1}^2((M_{\widetilde{p}}(\word{3}\word{j}\word{j})+M_{\widetilde{p}}(\word{j}\word{3}\word{j})+X_j(a)M_{\widetilde{p}}(\word{3}\word{j}))
\\&=2\sum_{j=1}^2((M_{\widetilde{p}}(\word{3}\word{j}\word{j})+M_{\widetilde{p}}(\word{3}\word{j}\shu\word{j}-2\cdot\word{3}\word{j}\word{j})+X_j(a)M_{\widetilde{p}}(\word{3}\word{j}))
\\&=2\sum_{j=1}^2(-M_{\widetilde{p}}(\word{3}\word{j}\word{j})+M_{\widetilde{p}}(\word{3}\word{j})\shu\word{j}+X_j(a)M_{\widetilde{p}}(\word{3}\word{j})).
\end{align*}
\end{example}

We are now ready to formulate  the main result of this section. Recall that $g\in\R[x_1,\dots,x_d]$ is a polynomial and $p:\R^d\to\R^{d+1}$ is the map introduced in \Cref{definition: p and M and phi and g}.

\begin{theorem}\label{theorem: vecchio risultato ipersuperficie}
Let $X:[a,b]\to \R^d$ be a reduced path and 
define  $\widetilde{g}(x)= g(x+X(a))-g(X(a))$. The following are equivalent:
\begin{enumerate}[label=(\arabic*),ref=(\arabic*)]
    \item\label{item1 thm carat ipersup} There exists a constant $c\in\R$ such that $\im(X)$ is contained in the hypersurface of equation $g(x_1,\dots,x_d)=c$.
\item\label{item2 thm carat ipersup} The signature of $X$ satisfies
\begin{equation}\label{eqCond2ThMAIN}
\left\langle\sigma(X),\sum_{j=1}^d\f_d\left(\frac{\partial \widetilde{g}}{\partial x_j}\right)\word{j} w\right\rangle=0 \mbox{ for every } w\in\TT(\R^d).
\end{equation}
\end{enumerate}
\end{theorem}
\begin{proof}Consider the map $p:\R^d\to\R^{d+1}$ as in  \Cref{definition: p and M and phi and g}. If $\im(X)$ is contained in  the hypersurface $g(x_1,\dots,x_d)=c$, then $\im(p\circ X)$ is contained in  the hyperplane $x_{d+1}=c$. By \cite[Proposition 6.3]{GalSan} or \cite[Lemma 3.1]{preiss}, this implies that \begin{align*}
&\langle\sigma(p\circ X),v\rangle=0 \mbox{ for every }v\in\TT(\R^{d+1})\mbox{ containing the letter }\word{d+1}.\\
&\Rightarrow \langle\sigma(p\circ X),v\rangle=0 \mbox{ for every }v\in\TT(\R^{d+1})\mbox{ that begins with the letter }\word{d+1}.\\
&\Rightarrow \langle\sigma(X),M_{\widetilde{p}}(v)\rangle=0 \mbox{ for every }v\in\TT(\R^{d+1})\mbox{ that begins with the letter }\word{d+1}
\end{align*}
by \Cref{theorem: Laura e Rosa}. The latter is equivalent to \ref{item2 thm carat ipersup} by \Cref{lemma: proprietà M}\eqref{item: le nostre equazioni sono M(ep)}.

The converse implication follows a similar argument. If \ref{item2 thm carat ipersup} holds, then   \Cref{lemma: proprietà M}\eqref{item: le nostre equazioni sono M(ep)} implies  that \begin{equation}
    \label{eq: zero se inizia con e}
\langle\sigma(X),M_{\widetilde{p}}((\word{d+1})w)\rangle=0 \mbox{ for every }w\in\TT(\R^{d}).
\end{equation} Now we want to prove the stronger statement that $\langle\sigma(X),M_{\widetilde{p}}(v)\rangle=0$ for every $v\in\TT(\R^{d+1})$  containing $\word{d+1}$. To this purpose we apply \Cref{corollary: any words is a shuffle combination of stuff with only one e at the beginning}: if $v\in\TT(\R^{d+1})$  contains $\word{d+1}$, then $M_{\widetilde{p}}(v)=M_{\widetilde{p}}((\word{d+1})q_1)\shu M_{\widetilde{p}}( r_1)+\dots+M_{\widetilde{p}}((\word{d+1})q_t)\shu M_{\widetilde{p}}(r_t)$, where none of the $q_i$ contains the letter $\word{d+1}$. Hence
\begin{align*}
\langle\sigma(X),M_{\widetilde{p}}(v)\rangle&=\sum_{i=1}^t\langle\sigma(X),M_{\widetilde{p}}((\word{d+1})q_i)\shu M_{\widetilde{p}}( r_i)\rangle
\\&=\sum_{i=1}^t\langle\sigma(X),M_{\widetilde{p}}((\word{d+1})q_i)\rangle\cdot\langle \sigma(X),M_{\widetilde{p}}( r_i)\rangle\mbox{ by \Cref{lem:shuffle identity}}
\\&=\sum_{i=1}^t 0\cdot\langle \sigma(X),M_{\widetilde{p}}( r_i)\rangle=0\mbox{ by \eqref{eq: zero se inizia con e}.}
\end{align*}
Hence we proved that $\langle\sigma(X),M_{\widetilde{p}}(v)\rangle=0$ for every $v\in\TT(\R^{d+1})$  containing $\word{d+1}$. By \Cref{theorem: Laura e Rosa}, this is equivalent to $
\langle\sigma(p\circ X),v\rangle=0$ for every $v\in\TT(\R^{d+1})$  containing $\word{d+1}$.
 Either by \cite[Proposition 6.3]{GalSan} or \cite[Lemma 3.1]{preiss}, this implies that $\im(p\circ X)$ is contained in the hyperplane $x_{d+1}=c$ for some constant $c\in\R$ and therefore $\im(X)$ is contained in the hypersurface $g(x_1,\dots,x_d)=c$.
\end{proof}

It is not difficult to generalize \Cref{theorem: vecchio risultato ipersuperficie} to any algebraic variety, not just hypersurfaces.
\begin{corollary}\label{coroll: caratterizzazione dei paths che stanno su una varietà}
Let $X:[a,b]\to \R^d$ be a reduced path. Let $g_1,\dots,g_t\in\R[x_1,\dots,x_d]$ and define  $\widetilde{g_i}(x)= g_i(x+X(a))-g_i(X(a))$. The following are equivalent:
\begin{enumerate}
    \item There exist constants $c_1,\dots,c_t\in\R$ such that $\im(X)$ is contained in the algebraic variety defined by the equations  $g_1=c_1,\dots,g_t=c_t$.
\item The signature of $X$ satisfies
\[\left\langle\sigma(X),\sum_{j=1}^d\f_d\left(\frac{\partial \widetilde{g_i}}{\partial x_j}\right)\word{j} w\right\rangle=0\mbox{ for every $i\in\{1,\dots,t\}$ and every } w\in\TT(\R^d).\]
\end{enumerate}
\end{corollary}

In \Cref{section: ODE} we will apply a slightly different version of \Cref{coroll: caratterizzazione dei paths che stanno su una varietà}: instead of allowing the hypersurface of equation $g_i$ to be translated by a constant $c_i\in\R$, we specify one of them by imposing that the hypersurface contains the initial point of the path.

\begin{corollary}\label{theorem: caratterizzazione dei paths che stanno su un'ipersuperficie algebrica}
Let $X:[a,b]\to \R^d$ be a reduced path and let $g_1,\dots,g_t\in\R[x_1,\dots,x_d]$. For each $i\in\{1,\dots,t\}$ define  $\widetilde{g_i}(x)= g_i(x+X(a))-g_i(X(a))$. The following are equivalent:
\begin{enumerate}
    \item $\im(X)$ is contained in the algebraic variety defined by the equations  $g_1=0,\dots,g_t=0$.
\item $g_1(X(a))=  \dots=g_t(X(a))=0$ and
\[\left\langle\sigma(X),\sum_{j=1}^d\f_d\left(\frac{\partial \widetilde{g_i}}{\partial x_j}\right)\word{j} w\right\rangle=0\mbox{ for every $i\in\{1,\dots,t\}$ and every } w\in\TT(\R^d).\]
\end{enumerate}
\end{corollary}

\begin{example}\label{example: hyperbolic paraboloid}
Let $X:[0,1]\rightarrow \R^3$ be a reduced path. Take $g=2x^2-y^2-z$, so  $\widetilde{g}=2x^2+4X_1(0)x-y^2-2X_2(0)y-z$. By \Cref{theorem: vecchio risultato ipersuperficie}, $\im(X)$ lies  on the hyperbolic paraboloid defined by $g(x,y,z)=c$ for some $c\in\R$ if and only if
\begin{equation}\label{eq: paraboloide}
\langle \sigma(X), (4\cdot {\word{1}\word{1}}-2\cdot {\word{2}\word{2}})w \rangle+4X_1(0)\langle \sigma(X), \word{1}w \rangle-2X_2(0)\langle \sigma(X), \word{2}w \rangle=\langle \sigma(X),\word{3}w\rangle
\end{equation}
for every $w\in \mathbb{T}(\mathbb{R}^3)$. If we apply instead  \Cref{theorem: caratterizzazione dei paths che stanno su un'ipersuperficie algebrica}, we can say that $\im(X)$ lies  on the hyperbolic paraboloid defined by $g(x,y,z)=0$ if and only if $g(X(0))=0$ and \eqref{eq: paraboloide} holds for every $w\in \mathbb{T}(\mathbb{R}^3)$.
\end{example}

\begin{remark}\label{remark: delle nostre equazioni non ne bastano un numero finito} One might wonder if the infinitely many conditions \eqref{eqCond2ThMAIN} of \Cref{theorem: vecchio risultato ipersuperficie} could be replaced by only finitely many. However, in general this is not possible: given $k\in\N$, we can always approximate the $k$-th signature of a path with the $k$-th signature of a piecewise linear path. More precisely, as pointed out at the beginning of the proof of \cite[Theorem 5.6]{AFS19}, for every smooth path $X:[a,b]\to\R^d$ and every $k\in\N$ there exists $m\in\N$ and a piecewise linear path $Y:[a,b]\to\R^d$ with $m$ steps such that $\sigma^{(k)}(Y)=\sigma^{(k)}(X)$. So if we fix an hypersurface of equation $g(x_1,\dots,x_d)=0$, say for instance a hypersphere, then there exist plenty of piecewise linear paths $Y$ that satisfy \eqref{eqCond2ThMAIN} for all words $w$ of length up to $k$, but do not lie on the hypershpere, because the hypersphere does not contain any segment.
\end{remark}

A characterization of paths lying on an algebraic hypersurface also appears in \cite{preiss}. The core difference between \Cref{theorem: vecchio risultato ipersuperficie} and \cite[Corollary 3.6]{preiss} is how they deal with translations. The latter aims to characterize the set
\[V_g=\{X:[a,b]\to\R^d\mbox{ path}\mid \im(X-X(a))\subseteq (g=0)\},\]
while \Cref{theorem: vecchio risultato ipersuperficie}  aims to characterize the set
\[W_g=\{X:[a,b]\to\R^d\mbox{ path}\mid \im(X)\subseteq (g=c)\mbox{ for some }c\in\R\}.\]
We compare them in two examples.

\begin{example}\label{example: confronto con Rosa - iperpiano}
Take $g=x_1$. According to \cite[Lemma 3.1]{preiss}, the set $V_{g}$ equals $$\{X:[0,1]\to\R^d\mbox{ path}\mid \langle\sigma(X),w\rangle=0\mbox{ for every word $w$ containing the letter \word{1}}\}.$$
On the other hand 
$\widetilde{g}=g$ and \Cref{theorem: vecchio risultato ipersuperficie} tells us that the set $W_{g}$ equals $$\{X:[0,1]\to\R^d\mbox{ path}\mid \langle\sigma(X),w\rangle=0\mbox{ for every word $w$ beginning with the letter \word{1}}\}.$$ Since the sets $V_g$ and $W_g$ are the same, our result characterizes the same set of paths by using less equations. This also improves \cite[Proposition 6.3]{GalSan}.
\end{example}

\begin{example}\label{example: confronto con Rosa - circle}
Let $g(x_1,x_2)=x_1^2-2x_1+x_2^2$. Let $I$ be the smallest vector subspace of $\TT(\R^2)$ containing the word $\word{1}\word{1}-\word{1}+\word{2}\word{2}$ and closed with respect to $\succ$. According to \cite[Example 3.8]{preiss}, the set $V_g$ is $$\{X:[0,1]\to\R^2\mbox{ path}\mid \langle\sigma(X),w\rangle=0\mbox{ for every }w\in I\}.$$
On the other hand 
$\widetilde{g}(x_1,x_2)=x_1^2+2X_1(0)x_1-2x_1+x_2^2+2X_2(0)x_2$ and \Cref{theorem: vecchio risultato ipersuperficie} tells us that $W_g$ is the set $$\{X:[0,1]\to\R^2\mbox{ path}\mid \langle\sigma(X),(\word{1}\word{1}+X_1(0)\word{1}-\word{1}+\word{2}\word{2}+X_2(0)\word{2})w\rangle=0\mbox{ for every } w\in\TT(\R^2)\}.$$ Unlike in \Cref{example: confronto con Rosa - iperpiano}, this time the sets $V_g$ and $W_g$ are not the same. 
However, when we restrict to paths starting at the origin, 
the two sets coincide. \Cref{theorem: vecchio risultato ipersuperficie} characterizes this set as
$$\{X:[0,1]\to\R^2\mbox{ path}\mid \langle\sigma(X),(\word{1}\word{1}-\word{1}+\word{2}\word{2})w\rangle=0\mbox{ for every } w\in\TT(\R^2)\},$$
which are less equations than the ones given by $I$.
\end{example}

One major difference is that $W_g$ is closed under taking the inverse path, while $V_g$ is not. The advantage of \cite[Corollary 3.6]{preiss} is that $V_g$ has nicer equations, involving only entries of the signature - in other words, $V_g$ is a \emph{path variety} as defined in \cite[Section 1]{preiss} - while equations of $W$ involve not only entries of $\sigma(X)$, but also entries of the vector $X(a)$. 
\begin{remark}
When we restrict our attention to paths starting at the origin, the subset $\{X\in V_g\mid X(a)=0\}$ coincides with $\{X\in W_g\mid X(a)=0\}$. For this subset of paths, the entries of $X(a)$ disappear and \Cref{theorem: vecchio risultato ipersuperficie} provides a characterization that uses less equations that the one provided by \cite[Corollary 3.6]{preiss}, as illustrated in \Cref{example: confronto con Rosa - circle}.
\end{remark}

\section{Signatures of holonomic paths}\label{section: lagrangian paths}

Legendrian paths -- or Legendrian \textit{knots}, as they are called in the specialized  literature -- arise as the integral curves of the \textit{contact distribution} in $\R^3$ with coordinates $x,y,p$. As illustrated in \cite[Section~2.4]{MR3760967}, this distribution is the kernel of the one--form $dy-pdx$, and it is the unique completely non--integrable  $\SL(2,\R)$--invariant distribution in $\R^3$. 
Because of complete non--integrability, Legendrian paths cannot be characterized as curves lying on a foliation of $\R^3$ by hypersurfaces, so we cannot use the same arguments as in \Cref{section: paths su una varietà}. In this section we take a slightly different approach and we characterize not only Legendrian, but also holonomic paths in terms of their signature tensors. This is the content of \Cref{theorem: characterization of holonomic paths}.

\begin{definition}\label{defLagCoord}\label{definition: legendrian}
A  path $X:[a,b]\to\R^3$  of class $\PW{1}$ is \emph{Legendrian} if 
$ \dot{X}_1(t)X_3(t)=\dot{X}_2(t)$ for every $t\in [a,b]$ for which $\dot{X}(t)$ exists.
\end{definition}

We start by providing a necessary condition.
\begin{lemma}\label{lemma: legendrian implies eq on signatures}
Let $X:[a,b]\to\R^3$ be a path of class $\PW{1}$. If  $X$ is  Legendrian, then
\begin{equation}\label{eqLegendrianita}
    \langle \sigma(X), \word{3}\word{1}w \rangle -\langle \sigma(X), \word{2}w \rangle+X_3(a)\langle \sigma(X), \word{1}w \rangle=0 \mbox{  for every } w \in \mathbb{T}(\mathbb{R}^3).
\end{equation}
\end{lemma}
\begin{proof}
Let $k\ge 2$ and let $ w =\word{i}_1\cdots \word{i}_{k-2}\in \mathbb{T}(\mathbb{R}^3)$, with the convention that if $k=2$ then $w = \varnothing$. 
By definition,
\begin{align*}
\langle \sigma(X), \word{31i}_1\cdots \word{i}_{k-2} \rangle=& \int_a^b \int_a^{t_k}\cdots \int_a^{t_2} \dot{X}_3(t_1)\dot{X}_{1}(t_2)\dot{X}_{i_1}(t_3)\cdots \dot{X}_{i_{k-2}}(t_k){\rm d}t_1 \cdots {\rm d}t_k \\
   =&\int_a^b \int_a^{t_k}\cdots \int_a^{t_3} \dot{X}_{1}(t_2)\dot{X}_{i_1}(t_3)\cdots \dot{X}_{i_{k-2}}(t_k) \int_a^{t_2} \dot{X}_3(t_1){\rm d}t_1 \cdots {\rm d}t_k\\
   =& \int_a^b \int_a^{t_k}\cdots \int_a^{t_3}
   X_3(t_2) \dot{X}_{1}(t_2)\dot{X}_{i_1}(t_3)\cdots \dot{X}_{i_{k-2}}(t_k)  {\rm d}t_2 \cdots {\rm d}t_k \\
   &- X_3(a)\int_a^b \int_a^{t_k}\cdots \int_a^{t_3}
  \dot{X}_{1}(t_2)\dot{X}_{i_1}(t_3)\cdots \dot{X}_{i_{k-2}}(t_k)  {\rm d}t_2 \cdots {\rm d}t_k.
 \end{align*} 
  Since $X$ is Legendrian, in the first summand we replace $ X_3(t_2) \dot{X}_{1}(t_2)$ by $\dot{X}_2(t_2) $. Moreover, we recognize that the second summand is $X_3(a)\langle \sigma(X), \word{1i}_1\cdots\word{i}_{k-2} \rangle$. Hence $\langle \sigma(X), \word{31i}_1\cdots\word{i}_{k-1} \rangle$ equals
 \begin{align*}
&\int_a^b \int_a^{t_k}\cdots \int_a^{t_3}
    \dot{X}_{2}(t_2)\dot{X}_{i_1}(t_3)\cdots \dot{X}_{i_{k-2}}(t_k)  {\rm d}t_2 \cdots {\rm d}t_k -X_3(a)\langle \sigma(X), \word{1i}_1\cdots\word{i}_{k-2} \rangle\\
&   = \langle \sigma(X), \word{2i}_1\cdots\word{i}_{k-2} \rangle-X_3(a)\langle \sigma(X), \word{1i}_1\cdots \word{i}_{k-2} \rangle.\qedhere
\end{align*}\end{proof}

By Lemma \ref{lemma: legendrian implies eq on signatures}
    condition \eqref{eqLegendrianita} is necessary for a path $X$ to be Legendrian: as a consequence of the main Theorem \ref{theorem: characterization of holonomic paths}, it will turn out to be sufficient as well, see Corollary \ref{theorem: characterization legendrian R3} below.

\begin{corollary}\label{corollario: legendrian implies eq on signatures}
Let $d\ge 3$ and let $X:[a,b]\to\R^d$ be a path of class $\PW{1}$. If there exist $p,q,s\in\{1,\dots,d\}$ such that the projection $\widetilde{X}(t)=(X_p(t),X_q(t),X_s(t))$ is  Legendrian, then $$ \langle \sigma(X), \word{s}\word{p}w \rangle -\langle \sigma(X), \word{q}w \rangle+X_s(a)\langle \sigma(X), \word{p}w \rangle=0 \mbox{  for every } w \in \mathbb{T}(\mathbb{R}^d).$$
\end{corollary}

Since the notion of Legendrian path only makes sense in $\R^3$, we want now to move to holonomic paths. Let $r\ge 1$ and $l\ge 0$ be integers and set $$d=1+r(l+1).$$

\begin{definition}\label{defRLholo}
For positive $l$, a path $X:[a,b]\rightarrow \mathbb{R}^{d}$   of class $\PW{1}$ is $(r, l)$-\emph{holonomic} if
\begin{equation*}
\label{eqHolRL}
\dot{X}_i(t)=\dot{X}_1(t)X_{i+r}(t), \hbox{ for all } i\in\{2,\dots,r l+1\}.
\end{equation*}
For $l=0$ we define every path of class $\PW{1}$ to be $(r, 0)$-holonomic.
\end{definition}
We observe that, for positive $l$, a path is $(r,l)$-holonomic if and only if,  for every $i\in\{2,\dots,r l+1\}$, the path  $[a,b]\to\R^3$ defined by $t\mapsto (X_1(t),X_i(t),X_{i+r}(t))$ is Legendrian.\par

In the geometric theory of PDEs, the canonical example of a holonomic path is the so--called $l$-th jet extension of a function, which  we define in the example below; see also  \cite[Section~4.1]{Geometry1}.  
\begin{example}\label{example: a path defined by derivatives of a function f is holonomic}
Let $k$ be a positive integer. If $f:[a,b]\to\R^r$ is a function of class $\PW{(k+l)}$, then we can define a $(r,l)$-holonomic path $X:[a,b]\to\R^d$ of class $\PW{k}$ by
\[X(t)=(t,f(t),f'(t),\dots,f^{(l)}(t)),
\]
known also as the $l$--order jet of $f$.
Indeed, if $i\in\{2,\dots,r l+1\}$, then there exists a unique $j\in\{0,\dots,l\}$ such that $i\in\{jr+2,\dots,(j+1)r+1\}$. Then $X_i(t)=f^{(j)}(t)_{i-jr-1}$ and $X_{i+r}(t)=f^{(j+1)}(t)_{i+r-(j+1)r-1}$. Hence \begin{align*}
\dot{X}_i(t)-\dot{X}_1(t)X_{i+r}(t)&=\frac{{\rm d}}{{\rm d}t}f^{(j)}(t)_{i-jr-1}-1\cdot f^{(j+1)}(t)_{i-jr-1}=0.
\end{align*}
The $(r,0)$--holonomic path $X(t)=(t,f(t))$ is nothing but the parametrized graph of $f$, known also as the zero--order jet of $f$.
\end{example}
 
\begin{remark}\label{remark: ridotto di un olonomo}\label{remark: rdotto di un legendriano è legendriano}
If $X$ is a Legendrian path, then the restriction of $X$ to any subinterval is Legendrian and the reduced path associated to $X$ is also Legendrian. In the same way, if $X$ is a $(r,l)$-holonomic path, then the reduced path associated to $X$ is also $(r,l)$-holonomic.
\end{remark}

Now we can state and prove our characterization of holonomic paths.
\begin{theorem}\label{theorem: characterization of holonomic paths}
Let $l\ge 0$ and $r\ge 1$ be integers and let $d=1+r(l+1)$. A  reduced path $X:[a,b]\to\R^{d}$ of class $\PW{1}$ is $(r, l)$-holonomic if and only if 
\begin{equation}\label{eq: for holonomic}
    \langle \sigma(X), { (\word{i}+\word{r})\word{1}}w \rangle -\langle \sigma(X),  \word{i}w \rangle+X_{i+r}(a)\langle \sigma(X),  \word{1}w \rangle=0
\end{equation}
for every every $w\in\TT(\R^d)$
and every integer $i$ such that $2\le i\le rl+1$
.
\end{theorem}
\begin{proof}
Assume that $X$ is $(r,l)$-holonomic. If $l=0$, then there is nothing to prove. If $l\ge 1$, then for  every $i\in\{2,\dots,rl+1\}$, the path $t\to (X_1(t),X_i(t),X_{i+r}(t))$ is Legendrian. By \Cref{corollario: legendrian implies eq on signatures}, $X$ satisfies \eqref{eq: for holonomic}.\par 

To prove the converse, we assume now  that $X$ satisfies \eqref{eq: for holonomic} and we introduce an auxiliary path $Y(t)=(Y_1(t),\ldots,Y_d(t))$ defined as follows: 
\begin{eqnarray*}
Y_1(t) &=& X_1(t)\, ,\\
Y_i(t) &=& X_i(a)+\int_a^tX_1'(q)Y_{i+r}(q)dq\mbox{ for every }i: 2\le i\le rl+1,\\
Y_i(t) &=& X_i(t)\mbox{ for every } i\in\{rl+2,\ldots, r(l+1)+1\}.
\end{eqnarray*}
Notice that if $l=0$ the above conditions reduce to  $Y=X$.\par 
Observe that $X(a)=Y(a)$ and that $Y$ is $(r, l)$-holonomic. By \Cref{corollario: legendrian implies eq on signatures}, $Y$ satisfies
\begin{equation}\label{eq: for holonomic Y}
    \langle \sigma(Y), { (\word{i}+\word{r})\word{1}}w \rangle -\langle \sigma(Y),  \word{i}w \rangle+Y_{i+r}(a)\langle \sigma(Y),  \word{1}w \rangle=0
\end{equation}
for every every $w\in\TT(\R^d)$
and every integer $i$ such that $2\le i\le rl+1$. We want to prove that $\sigma(X)=\sigma(Y)$
. We argue by induction on $l$. If $l=0$, then there is nothing to prove. 
Assume that $l\ge 1$ and let $\widetilde{X}:[a,b]\to\R^{d-2}$ be the projection of $X$ obtained by forgetting the components $X_2,X_3,\dots,X_{r+1}$. In other words, $\widetilde{X}(t)=(X_1(t),X_{r+2}(t), \dots,X_{r(l+1)+1}(t))$. In the same way take $\widetilde{Y}(t)=(Y_1(t),Y_{r+2}(t), \dots,Y_{r(l+1)+1}(t))$. Then $\widetilde{Y}$ is $(r,l-1)$-holonomic.\par  Let us prove that $\sigma(\widetilde{X})=\sigma(\widetilde{Y})$. If $l=1$, then this is true because  $\widetilde{X}=\widetilde{Y}$. Assume that $l\ge 2$. Now we prove that, since 
$X$ satisfies \eqref{eq: for holonomic} for every  $i\in\{2,\dots,rl+1\}$ and every $w\in\TT(\R^d)$, $\widetilde{X}$ satisfies \eqref{eq: for holonomic} for every  $j\in\{2,\dots,r(l-1)+1\}=\{2,\dots,d-r\}$ and every $w\in\TT(\R^{d-r})$. Indeed, if we call $\tau:\{1,2,\dots,d-r\}\to\{1,2,\dots,d\}$ the map $\tau(1)=1$ and $\tau(i)=i+r$ for every $i\in\{2,\dots,d-r\}$, and we take $j\in\{2,\dots,r(l-1)+1\}$,
\begin{align*}
\langle \sigma(\widetilde{X}),  \word{j}w \rangle&=\langle \sigma(X),  \tau(\word{j}w) \rangle=\langle \sigma(X),  (\word{j+r})\tau(w) \rangle\\
&=\langle \sigma(X),  (\word{j+2r})\word{1}\tau(w) \rangle+X_{j+2r}(a)\langle \sigma(X),  \word{1}\tau(w) \rangle\\
&=\langle \sigma(X),  \tau((\word{j+r})\word{1}w)) \rangle+\widetilde{X}_{j+r}(a)\langle \sigma(X),  \tau(\word{1}w) \rangle\\
&=\langle \sigma(\widetilde{X}),  (\word{j+r})\word{1}w \rangle+\widetilde{X}_{j+r}(a)\langle \sigma(\widetilde{X}),  \word{1}w \rangle.
\end{align*}
By induction hypothesis, $\sigma(\widetilde{X})=\sigma(\widetilde{Y})$, that is
\begin{equation}\label{equation: w senza lettere piccole}
\langle \sigma(X),w\rangle=\langle \sigma(Y),w\rangle\ \forall w \mbox{ not containing any letter in } \{\word{2},\dots,\word{r+1}\}.
\end{equation}In order to prove that $\sigma(X)=\sigma(Y)$, we take a word $\word{i}_1\dots\word{i}_k$ in the alphabet $\{\word{1},\dots,\word{d}\}$ and we show that $\langle \sigma(X),w\rangle=\langle \sigma(Y),w\rangle$ by induction on 
\[n=\#\{j\in\{1,\dots,k\}\mid \word{i}_j\in\{\word{2},\dots,\word{r+1}\}\}.
\]
If $n=0$, then $w$ does not contain any letter in $\{\word{2},\dots,\word{r+1}\}$ and the statement holds by \eqref{equation: w senza lettere piccole}. Now suppose thet $n\ge 1$ and assume  that 
$\langle\sigma(Y),\word{i}_1\dots\word{i}_k\rangle=\langle\sigma(X),\word{i}_1\dots\word{i}_k\rangle$ whenever at most $n-1$ of the letters belong to $\{\word{2},\dots,\word{r+1}\}$. Now we take a word $\word{i}_1\dots\word{i}_k$ in which $n$ of the letters belong to $\{\word{2},\dots,\word{r+1}\}$ and we set $w=\word{i}_2\dots\word{i}_k$. We distinguish two cases.
\begin{enumerate}[label=(\alph*)]
\item\label{bullet a} Assume that $\word{i}_1\in\{\word{2},\dots,\word{r+1}\}$. Then at most $n-1$ letters of $w$ belong to $\{\word{2},\dots,\word{r+1}\}$, so our induction hypothesis implies that 
$\langle\sigma(X),(\word{i}_1\word{+r})\word{1}w\rangle=\langle\sigma(Y),(\word{i}_1\word{+r})\word{1}w\rangle\mbox{ and }\langle\sigma(X),\word{1}w\rangle=\langle\sigma(Y),\word{1}w\rangle$.
Now we apply the hypothesis that \eqref{eq: for holonomic} holds for our choice of $w$, thus
\begin{align*}
\langle\sigma(X),\word{i}_1w\rangle&=\langle\sigma(X),(\word{i}_1\word{+r})\word{1}w\rangle+X_{i_1+r}(a)\langle\sigma(X),\word{1}w\rangle\\
&=\langle\sigma(Y),(\word{i}_1\word{+r})\word{1}w\rangle+X_{i_1+r}(a)\langle\sigma(Y),\word{1}w\rangle\\
&=\langle\sigma(Y),(\word{i}_1\word{+r})\word{1}w\rangle+Y_{i_1+r}(a)\langle\sigma(Y),\word{1}w\rangle\mbox{ because } Y(a)=X(a)\\
&=\langle\sigma(Y),\word{i}_1w\rangle\mbox{ because $Y$ satisfies } \eqref{eq: for holonomic Y}.
\end{align*}
\item\label{bullet b} Now we consider the general case, in which $\word{i}_1$ is any letter. Since $n\ge 1$, there exists $j\in\{1,\dots,k\}$ such that $\word{i}_j\in\{\word{2},\dots,\word{r+1}\}$. By \Cref{lemma: ridurre da una d+1 all'inizio ad una d+1 ovunque}, there exist words $q_1,\dots,q_t, r_1,\dots,r_t$ and coefficients $\lambda_1,\dots,\lambda_t\in\R$ such that $\word{i}_1\dots \word{i}_k=\lambda_1\word{i}_jq_1\shu r_1+\dots+\lambda_t\word{i}_jq_t\shu r_t$. \Cref{lem:shuffle identity} implies that  
\begin{align*}
\langle\sigma(X),\word{i}_1\dots \word{i}_k\rangle
&= \sum_{i=1}^t\lambda_i
\langle\sigma(X),\word{i}_j q_i\rangle\langle\sigma(X),r_i\rangle.
\end{align*}
But in the previous case we have already proven our statement for words which begins with the letter $\word{i}_j$, so  $\langle\sigma(X),\word{i}_jq_i\rangle=\langle\sigma(Y),\word{i}_jq_i\rangle$. Moreover, at most $n-1$ of the letters of  $r_i$ belong to $\{\word{2},\dots,\word{r+1}\}$, so our induction hypothesis implies that $\langle\sigma(X),r_i\rangle=\langle\sigma(Y),r_i\rangle$. Hence $$\langle \sigma(X),w\rangle=\sum_{i=1}^t\lambda_i\langle\sigma(X),\word{i}_jq_i\rangle\langle\sigma(X),r_i\rangle=\sum_{i=1}^t\lambda_i\langle\sigma(Y),\word{i}_jq_i\rangle\langle\sigma(Y),r_i\rangle=\langle\sigma(Y),w\rangle$$
by \Cref{lem:shuffle identity}.
\end{enumerate}
Even now that we have proven that $\sigma(X)=\sigma(Y)$, we cannot conclude that $X=Y$, because we cannot guarantee that $Y$ is reduced. So instead we consider $Y_{red}$, which has the same signature as $Y$ and is also $(r,l)$-holonomic by \Cref{remark: ridotto di un olonomo}. We have $\sigma(X)=\sigma(Y)=\sigma(Y_{red})$. Since both $X$ and $Y_{red}$ are reduced, \Cref{theorem: chen uniqueness} implies that $X=Y_{red}$, so we conclude that $X$ is $(r,l)$-holonomic.
\end{proof}

Because of the importance of Legendrian knots, we stress that  \Cref{theorem: characterization of holonomic paths}, in the case $l=r=1$, provides the following  inverse implication of Lemma \ref{lemma: legendrian implies eq on signatures}.

\begin{corollary}\label{theorem: characterization legendrian R3}
Let $X:[a,b]\to\R^3$ be a reduced path of class $\PW{1}$. Then $X$ is Legendrian if and only if \eqref{eqLegendrianita} are satisfied. 
\end{corollary}

\begin{remark}\label{remark: being legendrian is not translation invariant}
Notice that the set of equations of \Cref{theorem: characterization of holonomic paths} characterizing holonomic paths are not a variety in the sense of \cite[Section 1]{preiss}. This comes as no surprise since being holonomic is not translation-invariant.
\end{remark}

\section{Applications to ODE}\label{section: ODE}
The purpose of this section is to employ the results of \Cref{section: paths su una varietà} and \Cref{section: lagrangian paths} to establish necessary and sufficient conditions on the signature $\sigma(X)$ for a path $X$ to satisfy a systems of ODEs, up to tree--like excursions. 
Our framework   relies on the geometric theory of differential equations. We refer the interested reader to \cite[Chapter~4]{Geometry1}. 

In order to make the exposition clearer, this section is split in two parts: in  Section~\ref{subRicapitolazioneGeomODEs} we recall the fundamental definitions and sketch the main results of a geometric theory of ODEs based on jet spaces \cite[ Chapter~1]{SCLDEqMP}, whereas the theorem characterising solutions in terms of their signatures is formulated and proved in  Section~\ref{subRisultatoPrincipaleSezioneQuattro}. We would like to remark that early examples of applications of the signature framework in the context of ODE can be found in
\cite[Section~9]{preiss} and \cite[Slide~31]{RosaTALK2024}.

\subsection{A geometric framework for systems of ODEs}\label{subRicapitolazioneGeomODEs}
We will work with  systems of ordinary differential equations  of order $l$ in $r$ dependent variables:    as in \Cref{section: lagrangian paths}, set $$d=1+(l+1)r$$
and consider a system of $m$ ODEs 
\begin{equation}\label{equation: system of ODEs}
\begin{cases}
F_1(t,\x(t),\x'(t),\dots,\x^{(l)}(t))=0\\
\vdots\\
F_m(t,\x(t),\x'(t),\dots,\x^{(l)}(t))=0
\end{cases}
\end{equation}
of order $l$. The symbol $\x(t)$ appearing in~\eqref{equation: system of ODEs} denotes the unknown function $\x:\R\rightarrow \R^r$ of the  variable $t$.
The starting point of a geometric theory of differential equations based on jet spaces is to replace all the arguments of the functions  $F_1,\dots ,F_m$ by a string of $d$ new 
variables $x_1,\ldots,x_d$, that is, by a point of the space $\R^d$. As we shall see, it is convenient to regard 
system~\eqref{equation: system of ODEs} as a subset 
\begin{equation}\label{defODEgeom}
   \mathcal{E}=\{ (x_1,\dots,x_d)\in\R^d\mid F_i(x_1,\dots,x_d)=0\ \forall  i\in\{1,\ldots, m\}\}\subset\R^{d}\, ,
\end{equation}
which is the geometric counterpart of the system of ODEs. In the dedicated literature the space $\R^d$ is usually denoted by $J^l(1,r)$ and it is referred to as the space of $l$--jets of functions from $\R$ to $\R^r$, see for instance \cite[Chapter~3, Section~2]{SCLDEqMP}.
In order to easily recover the usual notion of a solution it is convenient to single out paths whose first entry equals the independent variables.
\begin{definition}
A path $X:[a,b]\to\R^d$ is called \emph{projectable} if $X_1(t)=t$ for every $t\in [a,b]$.
\end{definition}

\begin{theorem}\label{teorema:geomODE}
Let $k$, $l$ and $r$ be positive integers, let $d=1+r(l+1)$. Consider the system \eqref{equation: system of ODEs} and the associated set $\mathcal{E}$ defined as in \eqref{defODEgeom}. Let $X:[a,b]\to\R^{d}$ be a  path   of  class $\PW{k}$. If 
\begin{enumerate}
    \item\label{item 1 geomODE} $X$ is $(r,l)$--holonomic,
    \item\label{item 2 geomODE}  $X$ is projectable,
    \item\label{item 3 geomODE}  $\im(X)\subset \mathcal{E}$,
\end{enumerate}
then the function $f:\R\to\R^r$ defined by $f(t)=(X_2(t),\dots,X_{r+1}(t))$ 
is a   solution  of class $\PW{(k+l)}$ of \eqref{equation: system of ODEs}. 
Conversely, if $f:\R\to\R^r$ 
is a solution of \eqref{equation: system of ODEs} of class $\PW{(k+l)}$
, then its $l$--jet extension \begin{equation}\label{equation: l jet extension}
X(t)=(t,f(t),f'(t),\dots,f^{(l)}(t))
\end{equation}
is  of class $\PW{k}$ and it satisfies the three above conditions.
\end{theorem}
\begin{proof}
This is a standard result in the theory of PDEs based on jet spaces (see, for instance,   \cite[Section~4.1]{Geometry1}): for the convenience of the reader we adapt the proof to the particular case of ordinary differential equations. If we assume that a path $X:[a,b]\to\R^{d}$ satisfies the three conditions, then
\begin{equation}\label{eqJetExtNonOlonomo}
    X(t)=(t,\x_0(t), \x_1(t),\x_2(t),\ldots, \x_l(t))\, ,
\end{equation}
where
\begin{equation}\label{eqDefICSpiccoliConI}
    \x_i(t)=(x_{i,1}(t),x_{i,2}(t),\ldots,x_{i,r}(t))\mbox{ for every } i\in\{0,\ldots, l\}.
\end{equation}
Condition \eqref{item 1 geomODE} gives
\begin{equation}\label{eqCondHolLRdimTeo}
    \dot{X}_a(t)=\dot{X}_1(t)X_{a+r}(t)=X_{a+r}(t)\mbox{ for every } a\in\{2,\dots,r l+1\}.
\end{equation}
Switching to notation~\eqref{eqDefICSpiccoliConI}, conditions~\eqref{eqCondHolLRdimTeo} reads
\begin{equation*}
    x_{i+1, j}(t)=x'_{i,j}(t) \mbox{ for every } i\in\{0,\ldots, l-1\} \mbox{ and every } j\in\{1,\ldots, r\}.
\end{equation*}
In other words  \begin{equation}\label{eqJetExt}
X(t)=(t,\x_0(t),\x_0'(t),\dots,\x_0^{(l)}(t)).
\end{equation} 
Finally, condition \eqref{item 3 geomODE} means that~\eqref{equation: system of ODEs} is satisfied by all $t\in [a,b]$. It remains only to observe that, since $X$ is   of  class $\PW{k}$, so are its last $r$ entries $\x_0^{(l)}(t)$, and then $\x_0(t)$ must be of class $\PW{(k+l)}$. We conclude by taking $f=\x_0$.\par 
Conversely, if $f(t)$ is a solution of class $\PW{(k+l)}$, then the path $X:[a,b]\to\R^{d}$ defined as in \eqref{equation: l jet extension} is of class $\PW{k}$, since  since we took derivatives of order up to $l$ as in \Cref{example: a path defined by derivatives of a function f is holonomic}. By construction, $X(t)$  is $(r,l)$--holonomic (condition \eqref{item 1 geomODE}) and projectable (condition \eqref{item 2 geomODE}) and since~\eqref{equation: system of ODEs} is satisfied by all $t\in[a,b]$, condition \eqref{item 3 geomODE} is met as well.
 \end{proof}
In short, Theorem~\ref{teorema:geomODE} says that $\x$ is a solution of $\mathcal{E}$ if and only if its $l$--jet extension $X$ is contained into $\mathcal{E}$: in this case, since   the correspondence between $f$ and $X$ is one--to--one (see Example \ref{example: a path defined by derivatives of a function f is holonomic}), it is not misleading  to call \emph{solution} both $f$ and its jet extension $X$.\par

In most physical problems, a Cauchy problem comes with a set of initial conditions. We want to formulate them geometrically as well.

\begin{example}\label{exSimplCP}\label{example: easy Cauchy problem}
Let $r=l=1$. Let $G:\R^2\to\R$ be a polynomial function and consider the Cauchy problem
\begin{equation*}\label{eqCauchyProb}
\left\{\begin{array}{rcl}
    x'(t) &=& G(t,x(t))\, ,\\
    x(0) &=& c.
    \end{array} \right.
\end{equation*}
In order to rephrase it, define
    $F(x_1,x_2,x_3)=x_3-G(x_1,x_2)\, .$
Thus we get the surface $\mathcal{E}=\{(x_1,x_2,x_3)\in\R^3\mid F(x_1,x_2,x_3)=0\}\subset\R^{1+r(l+1)}=\R^3$. The initial condition $x(0)=c$ singles out another subset
\begin{equation*}\label{eqDefSigmaEx}
    \Sigma=\{(x_1,x_2,x_3)\in\R^3\mid x_2=c\}.
\end{equation*}
If we consider the path $X:[0,1]\to\R^3$ defined by $X(t)=(t,x(t),x'(t))$, then the initial condition $x(0)=c$ can be recast  as $X(0)\in\Sigma$.
\end{example} 
 
 In \Cref{example: easy Cauchy problem}, we are requiring the point $X(0)=(0,x(0),x'(0))$ to belong to  the affine plane  
of equation $x_2=c$. We want to give a more general definition (which is still a particular case of the most general one, see \cite[Section~2.5]{Geometry1}
).\par 
 \begin{definition}\label{defInitSurf} In the context of a Cauchy problem, the affine subspace of the form
 \begin{equation}\label{defSIGMAgeom}
\Sigma=\{ (x_1,\dots,x_{d})\in\R^{d}\mid x_a=c_a\mbox{ for every } a\in\{2,3,\ldots, d-r\}\}
\end{equation}
given by the initial conditions is called 
an \textit{initial subset} or \emph{initial data}.
 \end{definition}

Let $a,b\in\R$, such that $a<0<b$, and let $k\geq 0$ be an integer.   We show now how to interpret  the pair  $(\mathcal{E},\Sigma)$
  as a  Cauchy problem.

\begin{definition}\label{defSolCPgeo}
 A   \emph{solution} of class $\PW{(k+l)}$ of a Cauchy problem $(\mathcal{E},\Sigma)$ is a   path $X:[a,b]\to \R^{d}$   of class $\PW{(k+l)}$, such that:
\begin{enumerate}
    \item $X$ is $(r,l)$--holonomic,
    \item\label{item 2 defSolCPgeo} $X$ is projectable,
    \item $\im(X)\subset \mathcal{E}$,
    \item\label{item: condizione iniziale sigma} $X(0)\in\Sigma$.
\end{enumerate}
If the condition \eqref{item 2 defSolCPgeo} above is not met, then $X$ is  a \emph{generalized solution}.
\end{definition}

\begin{figure}[h]
  \begin{subfigure}[b]{0.3\textwidth}
\includegraphics[width=\linewidth]{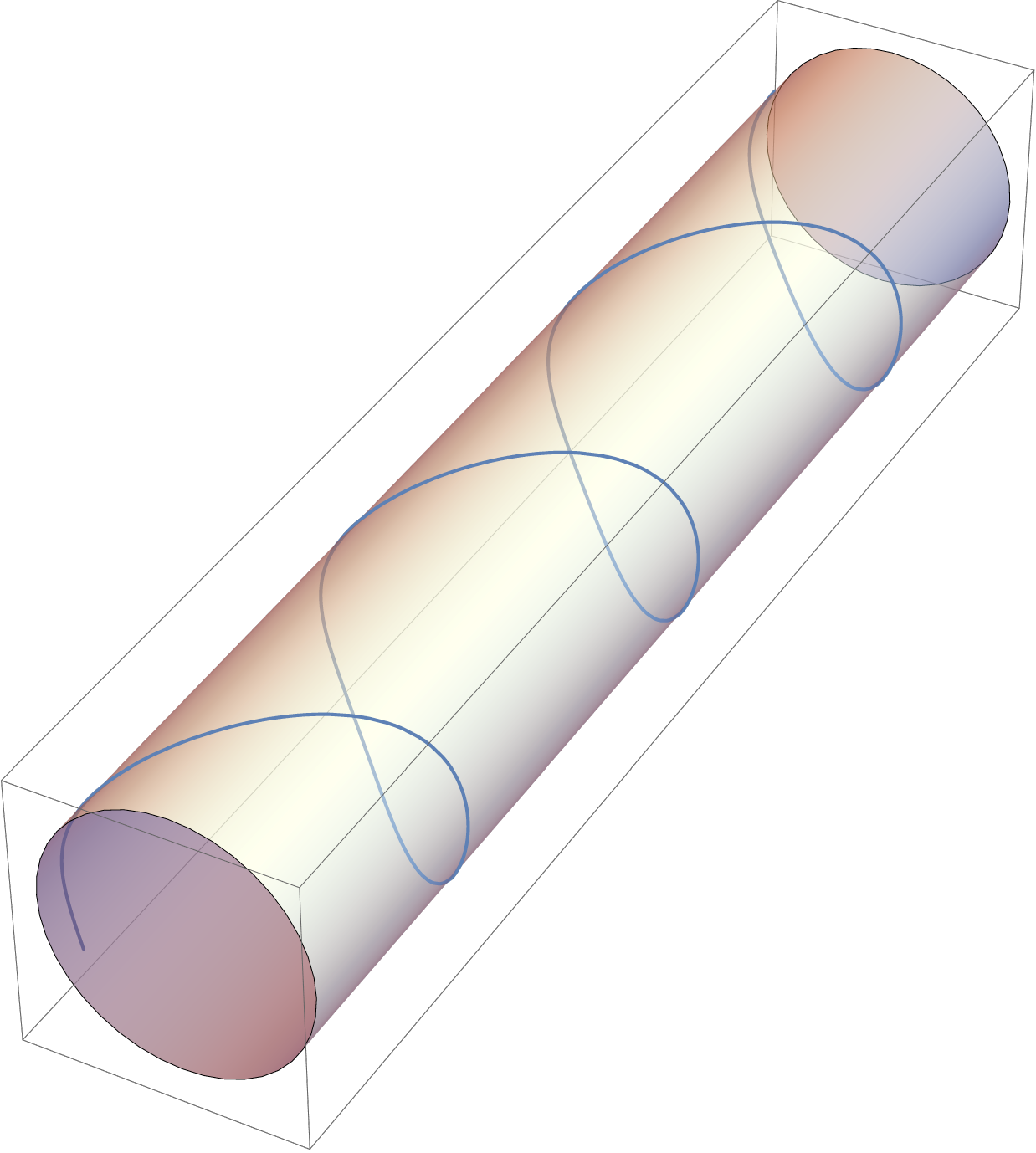}
    \caption{A Legendrian path on the surface of a cylinder.}
   \label{pic: Leg path on cylinder}
  \end{subfigure}
  \hfill
  \begin{subfigure}[b]{0.4\textwidth}
    \includegraphics[width=\linewidth]{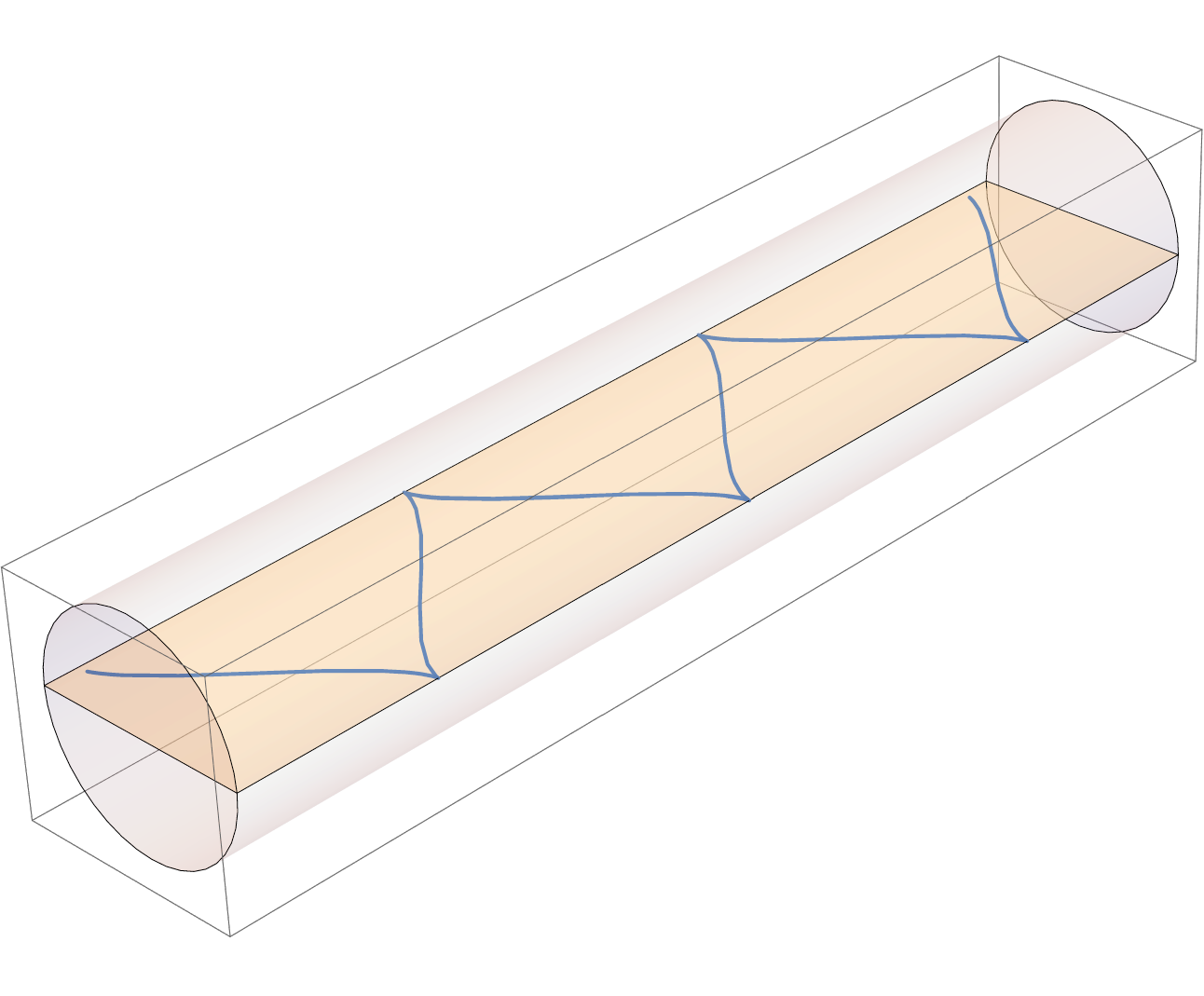}
    \caption{Projection of the Legendrian path of \Cref{pic: Leg path on cylinder}.}
    \label{pic: proj Leg path}
  \end{subfigure}
  \caption{\Cref{ex: path su cilindro}.}
\end{figure}

We give an example of a generalized solution of an ODE.
\begin{example}\label{ex: path su cilindro} 
Take $l=r=1$ and  $[a,b]=[-3\pi,3\pi]$,  and consider the ODE $x'(t)^2+t^2=1$.  The set $\mathcal{E}$ associated with it is the cylinder $\E=\{(x_1,x_2,x_3)\in\R^3\mid x_3^2+x_1^2=1\}$. If we use cylindrical coordinates, up to reparametrization, a path $X:[-3\pi,3\pi]\to\E$ is of the form
\begin{equation*}
    X(t)=(\cos(t),h(t),\sin(t))\, ,
\end{equation*}
for some 
function $h:[0,1]\to \R$. If $X$ is a Legendrian path, 
 then $h'(t)=-\sin^2(t)$. Easy integration leads to the generic Legendrian path  in $\mathcal{E}$:
\begin{equation}\label{eqCurvaSulCilindro}
    X(t)=(\cos(t),\tfrac{1}{2}(\sin(t)\cos(t)-t)+c,\sin(t))\mbox{ for some }c\in\R.
\end{equation}
 \Cref{pic: Leg path on cylinder} shows a Legendrian  $X$ on the surface of the cylinder $\mathcal{E}$. A brute--force projection of $X$ on the $(x_1,x_2)$--plane gives a path  in $\R^2$ depicted in \Cref{pic: proj Leg path},
which clearly is not the graph of any function.  By solving~\eqref{eqCurvaSulCilindro} in the $t$ variable, we find the paths 
\begin{equation}\label{eqSolExplicit}
    (t,\phi_{\pm}(t),\phi_{\pm}'(t))\, ,
\end{equation}
where
\begin{equation*}
    \phi_{\pm}(t)=\pm \tfrac{1}{2}(t\sqrt{1-t^2}\mp\arccos(t))+c.
\end{equation*}
Now it is not hard to see that, by concatenating six paths of the form~\eqref{eqSolExplicit}, with appropriate choices of $c$, one gets the above projection of $X(t)$. 
\end{example}

Much as \Cref{teorema:geomODE} showed that the usual notion of a solution $\x$ of a system of ODEs translates into $X$ being a path contained in $\mathcal{E}$,   Theorem~\ref{theorem: olonomo se e solo se è dato dalle derivate di f} below shows that     Definition~\ref{defSolCPgeo} is merely a    geometrical reformulation of the usual notion of a solution of Cauchy problem. Let $k$, $l$ and $r$ be positive integers and let $d=1+(l+1)r$. Let $a,b\in\R$ such that $a\le 0<b$. If $l\ge 2$, then we also require that $a<0$.

\begin{theorem}\label{thTeoremaFalso}\label{theorem: olonomo se e solo se è dato dalle derivate di f} 
Let $F_1,\dots,F_m\in\R[x_1,\ldots,x_d]$ be polynomials,    $\vv_0,\dots,\vv_l\in\R^r$, and consider the Cauchy problem
\begin{equation}\label{eqCPmolto_generale}
   \left\{\begin{array}{ll}
      F_i(t,\x(t), \x'(t),\ldots, \x^{(l)}(t))=0  &\mbox{ for every }i\in\{1,\dots,m\}  \\
\x^{(j)}(0)=\vv_j  & \mbox{ for every }j\in\{1,\dots,l-1\}.
   \end{array}\right.
\end{equation}
Then formula \eqref{equation: l jet extension} establishes a one--to--one correspondence between   solutions $f:[a,b]\to\R^r$ of class $\PW{(l+k)}([a,b])$ of \eqref{eqCPmolto_generale} and the solutions $X:[a,b]\to\R^{d}$  of class $\PW{k}([a,b])$ of the Cauchy problem $(\mathcal{E},\Sigma)$ in the geometric sense of Definition~\ref{defSolCPgeo}, where $\mathcal{E}$ and $\Sigma$ are given by~\eqref{defODEgeom} and~\eqref{defSIGMAgeom}, respectively.
\end{theorem}
\begin{proof}
A straightforward completion of the proof of \Cref{teorema:geomODE}. Indeed, the path $X(t)$ given by~\eqref{equation: l jet extension} satisfies the condition $X(0)\in\Sigma$, with $\Sigma$ given by~\eqref{defSIGMAgeom}, if and only if the $i$--th derivative $f^{(i)}(t)$ of $f(t)$ takes the value $\vv_i$ at the time $t=0$.  In other words, condition \eqref{item: condizione iniziale sigma} for $X(t)$ to be a solution of the Cauchy problem $(\mathcal{E},\Sigma)$ in  the geometric sense of Definition~\ref{defSolCPgeo} is equivalent to the fact that   $f(t)$ fulfills  the initial conditions of~\eqref{eqCPmolto_generale}.
\end{proof}
\begin{remark}\label{remark: holonomic is made by derivatives of f for generalized solutions}
If the projectability condition is dropped in Definition~\ref{defSolCPgeo},  then Theorem~\ref{thTeoremaFalso} is still valid because   a generalized solution can be broken into several honest solutions, leaving out a discrete set of singular points. In the setting of Example~\ref{ex: path su cilindro}, honest  solutions are   graphs of functions $y(t)$ and the aforementioned singular  points  corresponds to the   points at which the contact distribution is tangent to the surface $\mathcal{E}$.
\end{remark}

\subsection{Characterization of solutions in terms of signature}\label{subRisultatoPrincipaleSezioneQuattro}

Thanks to Theorem~\ref{thTeoremaFalso}, we can now work entirely within a geometric framework which allows to characterize a solution of a Cauchy problem in terms of signature tensors. Let $k, r,l$ be positive integers   and set  $d=1+r(l+1)$. Let $a,b\in\R$, such that $a\le 0<b$. If $l\ge 2$, then we also require that $a<0$.
\begin{theorem}\label{corollary: risultato generale di caratterizzazioni delle soluzioni di un'ODE}
Let $F_1,\dots,F_m\in\R[x_1,\ldots,x_d]$ be polynomials and $\vv_0,\dots,\vv_{l-1}\in\R^r$ be vectors.   
Set  $$\Sigma=\{(x_1,\dots,x_d)\in\R^d\mid (x_1,\dots,x_{d-r})=(0,\vv_0,\dots,\vv_{l-1})\}$$ 
and consider the 
Cauchy problem
\begin{equation}\label{equation: general cauchy problem}
\begin{cases}
F_1(t,\x(t),\x'(t),\dots,\x^{(l)}(t))=0\\
\vdots\\
F_m(t,\x(t),\x'(t),\dots,\x^{(l)}(t))=0\\
\x(0)=\vv_0 \\
\x'(0)=\vv_1 \\\vdots\\
\x^{(l-1)}(0)=\vv_{l-1}.
\end{cases}\end{equation}
If a   function $f:[a,b]\to\R^r$  of class  $\PW{(k+l)}$ is a solution of \eqref{equation: general cauchy problem}, then the path $X:[a,b]\to\R^
d$ defined by $X(t)=(t,f(t),f'(t),\dots,f^{(l)}(t))$  is reduced and
\begin{enumerate}
\item\label{item: main theorem, condizione iniziale} $X(0)\in\Sigma $,
\item\label{item: main theorem, olonomia} $\langle \sigma(X),  (\word{i}+\word{r})\word{1}w \rangle -\langle \sigma(X), \word{i}w \rangle+X_{i+r}(a)\langle \sigma(X), \word{1}w \rangle=0$ for every $i\in\{2,\dots,rl+1\}$ and for every  $w\in \mathbb{T}(\mathbb{R}^d
)$,
\item\label{item: main theorem, stare su tutti i polinomi} 
$\left\langle\sigma(X),\sum\limits_{j=1}^
d\f_d
\left(\frac{\partial \widetilde{F}_i}{\partial x_j}\right)\word{j} w\right\rangle=0$ for every $i\in\{1,\dots,m\}$ and for every $w\in\TT(\R^d
)$, where $\f_
d$ is defined in \Cref{definition: p and M and phi and g} and $\widetilde{F}_i(p)=F_i(p+X(a))-F_i(X(a))$ for every $p\in\R^d$.
\end{enumerate}

Conversely, if a   reduced  path $X:[a,b]\to\R^d$   of class     of class  $\PW{k}$ satisfies the three conditions, then the function $f(t)= (X_2(t),\dots,X_{r+1}(t))$ is a   generalized solution of \eqref{equation: general cauchy problem}   of class   $\PW{(k+l)}$.
\begin{proof}
Suppose that $f$ is a solution of \eqref{equation: general cauchy problem}. The initial conditions $\x^{(i)}(0)=\vv_i$ imply that $f^{(i)}(0)=\vv_i$ for every $i\in\{0,\dots,l-1\}$ and therefore $X(0)=(0,\vv_0,\vv_1,\dots,\vv_{l-1},f^{(l)}(0))\in\Sigma$, so condition \eqref{item: main theorem, condizione iniziale} holds. Second, since $X$ is $(r,l)$-holonomic (see \Cref{example: a path defined by derivatives of a function f is holonomic}),  condition \eqref{item: main theorem, olonomia} holds by \Cref{theorem: characterization of holonomic paths}. Moreover the differential equations $F_j(t,f(t),f'(t),\dots,f^{(l)}(t))=0$ imply that $F_j(X(t))=0$ for all $t\in[a,b]$ and for every $j\in\{1,\dots,m\}$. Hence condition \eqref{item: main theorem, stare su tutti i polinomi} holds by \Cref{coroll: caratterizzazione dei paths che stanno su una varietà}.

    Conversely, assume that a path  $X:[a,b]\to\R^d$ of class $\PW{k}$ satisfies conditions \eqref{item: main theorem, condizione iniziale}, \eqref{item: main theorem, olonomia} and \eqref{item: main theorem, stare su tutti i polinomi}. Then $X(0)\in\Sigma$, $X$ is $(r,l)$-holonomic by \Cref{theorem: characterization of holonomic paths} and $F_i(X(t))=0$ for all $t\in [a,b]$ and  for every $i\in\{1,\dots,m\}$ by \Cref{coroll: caratterizzazione dei paths che stanno su una varietà}. If we define $f:[a,b]\to\R^r$ by $f(t)=(X_2(t),\dots,X_{r+1}(t))$, then    \Cref{theorem: olonomo se e solo se è dato dalle derivate di f}, implies that  $f$ is a generalized solution (see  \Cref{remark: holonomic is made by derivatives of f for generalized solutions}) of the Cauchy problem \eqref{equation: general cauchy problem} of class $\PW{(k+l)}$.
\end{proof}
\end{theorem}

\begin{remark}\label{remark: le soluzioni di ode non sono translation invariant quindi le condizioni non sono rosee}
In general, none of the three conditions of \Cref{corollary: risultato generale di caratterizzazioni delle soluzioni di un'ODE} defines a path variety in the sense of \cite{preiss}, see also \Cref{remark: being legendrian is not translation invariant}. This does not come as a surprise, since the set of solutions of a Cauchy problem does not need to be translation--invariant.
\end{remark}

\begin{example}\label{proposition: caratterizzazione delle ODE in R3}
Let $l=r=1$, so that $d=3$. Let  $F=(x_1-1)^2+x_2^2+ x_3^2-1\in\R[x_1,x_2,x_3]$ be the sphere of radius 1 centered at $(1,0,0)$. We apply \Cref{corollary: risultato generale di caratterizzazioni delle soluzioni di un'ODE} to characterize solutions in $\PW{1}([0,1])$ of the Cauchy problem
\begin{equation*}\label{eq: ivp0}
\begin{cases}
x(t)^2+x'(t)^2=2t-t^2\\
x(0)=0
\end{cases}
\end{equation*}
in terms of signatures of the path $X:[0,1]\rightarrow \mathbb{R}^3$ defined by $X(t)=(t,x(t),x'(t))$. 
Condition \eqref{item: main theorem, condizione iniziale} of \Cref{corollary: risultato generale di caratterizzazioni delle soluzioni di un'ODE} says that $X(0)=(0,0,x'(0))$. Condition \eqref{item: main theorem, olonomia} of being Legendrian translates into 
\[\langle \sigma(X),  \word{31}w \rangle -\langle \sigma(X), \word{2}w \rangle+X_{3}(0)\langle \sigma(X), \word{1}w \rangle=0\mbox{ for every }w\in \mathbb{T}(\mathbb{R}^3
).\]
For condition \eqref{item: main theorem, stare su tutti i polinomi}, first observe that $\widetilde{F}(x_1,x_2,x_3)=x_1^2+x_2^2+x_3^2+2x'(0)x_3 $ and its derivatives are  $\nabla \widetilde{F}=(2x_1 ,\, 2x_2,\, 2(x_3+x'(0)))$. Hence condition \eqref{item: main theorem, stare su tutti i polinomi} translates to the signature of $X$ as
\[\langle \sigma(X),  (\word{11}+\word{22}+\word{33})w \rangle=x'(0)\langle \sigma(X), \word{3}w \rangle\mbox{ for every }w\in \mathbb{T}(\mathbb{R}^3
).\]  
\end{example}

\begin{remark}
In the second part of the proof of \Cref{corollary: risultato generale di caratterizzazioni delle soluzioni di un'ODE}, if one assumes $X_1(t)=t$, then in view of $\dot{X}_1(t)X_{i+r}(t)=\dot{X}_i(t)$, one obtains a solution $f(t)=(X_2(t),\dots,X_{r+1}(t))$.
\end{remark}

\section{Integral curves of vector fields}
\label{section: esempi e applicazioni}
A remarkable case of first--order ODEs that fits particularly well the framework of signatures is the equation describing the integral curves of a vector field $\Phi$ on $\R^r$, that is a global smooth section $   \Phi\in\Gamma(\R^r,T(\R^r))$ of the tangent bundle. As illustrated in \cite[Section~3]{MR1202431},  $\Gamma(\R^r,T(\R^r))$ is a free $\CC^\infty(\R^r)$--module of dimension $r$, generated by the partial derivatives $\tfrac{\partial}{\partial x_1},\dots,\tfrac{\partial}{\partial x_r}$. For any $p\in\R^r$, the value of the section $\Phi$ at the point $p$ is denoted by
$$
\Phi|_{p}\in\  T_{p}\R^r\simeq\R^r\, .
$$
\begin{definition}
A curve $\x:  [a,b]\to\R^r$ is called an \emph{integral curve} of the vector field $\Phi$ in $[a,b]$ if 
\begin{equation*}\label{eqFlowVF1}
\x'(t)=\Phi|_{\x(t)}\mbox{ for every } t\in [a,b],
\end{equation*}
that is, the tangent vector $\x'(t)$ to the curve ${\x}$ at the time $t$ equals the value of the section $\Phi$ at the point ${\x}(t)$.
\end{definition}
\begin{remark}\label{rem: integral curve as system of ode}
A vector field $\Phi\in \Gamma(\R^r,T(\R^r))$ can be written as a linear combination
\begin{equation*}\label{eqCampVettCoord}
    \Phi=\sum_{i=1}^r\Phi_i\frac{\partial }{\partial x_i}
\end{equation*}
of the basis elements of $\Gamma(\R^r,T(\R^r))$. Then the definition of an integral curve is a system of $r$ first--order ODEs
\[
\begin{cases}
\x'_1(t)=\Phi_1(\x(t))\\
\vdots\\
\x'_r(t)=\Phi_r(\x(t)).
\end{cases}
\]
In our notation \eqref{equation: system of ODEs} we have $l=1$, $d=1+2r$ and  
\begin{equation}\label{eqDeEffeI}
F_i(t,\x,\y)=\Phi_i(\x)-y_i.
\end{equation}
\end{remark}

 \subsection{Linear vector fields}
When all the smooth functions $\Phi_i$ appearing in \Cref{rem: integral curve as system of ode} 
are linear, each function
 \begin{equation*}
     F_i:\R^d=\R\oplus\R^r\oplus\R^r\to\R
 \end{equation*}
 appearing in \eqref{eqDeEffeI} does not depend upon the first constituent  $\R$, is linear on the second constituent $\R^r$, and $F_i(0,0,(y_1,\dots,y_r))=-y_i$
 . 
In such a case $\Phi$ is   an endomorphsim of $\R^r$ and it is called a \textit{linear vector field}. Then, in the standard basis of $\R^r$, 
 we identify $\Phi$ with a matrix $A=(a_{i,j})\in\R^{r\times r}$. 
 Accordingly, the system of ODEs defining the integral curve in \Cref{rem: integral curve as system of ode} becomes 
 \begin{equation*}\label{eqFlowVF1-2LIN}
\x'_i(t)=\sum_{j=1}^ra_{i,j}x_j(t)\mbox{ for every } i\in\{1,\ldots,r\}.
\end{equation*}
In matrix form, $ \x'(t)=A\cdot \x(t)$ for all $t$, and \eqref{eqDeEffeI} reads
\begin{equation*}\label{eqDeEffeILIN}
F_i(t,\x,\y)=\sum_{j=1}^ra_{i,j}x_j(t)-y_i\mbox{ for every } i\in\{1,\ldots,r\}.
\end{equation*}
As an example, when we  take  $A$ to be a skew--symmetric matrix, the integral curves  will travel along hyperspheres on $\R^r$ and then form the variety of paths lying on a hypersphere centered at the origin.\par 

\begin{example}
Let $A\in \mathbb{R}^{r\times r}$ be skew-symmetric, and let  $\x:[0,1]\rightarrow \R^r$ be an integral curve  of the linear vector field determined by $A$. We want to
prove that $\mathrm{im}(\x) $ lies on the hypersphere $x_1^2+\cdots +x_r^2=\|\x(0)\|^2$. Set $g(x_1,\dots,x_r)=x_1^2+\cdots +x_r^2-\|\x(0)\|^2$. By \Cref{theorem: caratterizzazione dei paths che stanno su un'ipersuperficie algebrica}, $\mathrm{im}(X)$ lies on the hypersphere if and only if $$\sum_{i=1}^r \langle \sigma(X),\word{i}w \rangle+X_i(0)\langle\sigma(X),w  \rangle =0, \hbox{ for all words }w.$$
The above equality can be proved with a direct computation using that $0=\x(t)^TA\x(t)=\x(t)^T\x'(t)=\x(t)\cdot \x'(t)$.
The above conditions are necessary conditions for a path  $\x(t)$ to be an integral curve of the linear vector field associated with $A$.
To obtain sufficient conditions, one has to rely on \Cref{corollary: risultato generale di caratterizzazioni delle soluzioni di un'ODE} and, in particular, work on the extension $X(t)$ of $\x(t)$ to $\R^d$. Indeed, 
from a geometric point of view, the space of skew--symmetric matrices of size $r$ is the Lie algebra $\mathfrak{so}_r(\R)$ of the Lie group of orthogonal transformations $\mathrm{SO}_r(\R)$ of $\R^r$. Integral curves, in this case, are described by the exponential map $A_t=\exp(tA)$ and they turn out to be rotations. All rotations lies on a hypersphere, but not all paths lying on a hypersphere are rotations.
\end{example}

According to \Cref{rem: integral curve as system of ode}, 
a  curve is integral if it is a solution of a certain system of ODE. When such equations are polynomial,  \Cref{corollary: risultato generale di caratterizzazioni delle soluzioni di un'ODE} characterizes integral curves according to the signature of a suitable path. We give the precise statement for the case of linear vector fields.

\begin{corollary}\label{corollary: linear vf}
Let $p\in\R^r$. Let $A=(a_{i,j})\in \R^{r\times r}$ and set $d=1+2r$. 
If $f:[a,b]\to\R^r$, with $a\le 0<b$, is the integral curve of the linear vector  field   defined by $A$ fulfilling  $f(0)=p$, then
the path $X:[a,b]\to\R^d$ defined by  $X(t)=(t,f(t),f'(t))$ satisfies 
\begin{align}\label{eq:punto iniziale per linear integral vector curve}
X(0)=(0,p,A\cdot p) \hbox{ and }
\sum_{j=1}^ra_{i,j}\langle\sigma(X),(\word{1}+\word{j})w \rangle=\langle \sigma(X),(\word{1}+\word{r}+\word{i})w  \rangle 
\end{align}
for every $i\in\{1,\dots,r\}$ and every word $w\in\TT(\R^d)$. Conversely, if a  reduced path $X:[a,b]\to \R^d$ is $(1,r)$-holonomic and satisfies \eqref{eq:punto iniziale per linear integral vector curve},  
then the curve $f:\R\to\R^r$ defined by $f(t)=(X_2(t)\dots,X_{r+1}(t))$ is integral with respect to the linear vector field defined by $A$.
\end{corollary}
\begin{proof}
For $i\in\{1,\dots,r\}$, define the linear polynomials
\begin{equation*}\label{eqFIcampoLineare}
F_i(t,x_1,\dots,x_r,y_1\dots,y_r)=\sum_{j=1}^r a_{i,j}x_j-y_i.\end{equation*}
Let $f$ be a solution of the system of ODEs
$$
\begin{cases}
F_1(t,\x(t),{\x'}(t))=0\\
\vdots \\
F_r(t,\x(t),{\x'}(t))=0
\end{cases}
$$
and let $X(t)=(t,f(t),f'(t))$. A direct computation shows that $\widetilde{F}_i=F_i$ for all $i\in \{1,\dots,r \}$.
Therefore 
\begin{align*}
\frac{\partial \widetilde{F_i}}{\partial t}=0\mbox{, }\frac{\partial \widetilde{F_i}}{\partial \x_j}=a_{i,j}\mbox{ and } \frac{\partial \widetilde{F_i}}{\partial \y_j}=-\delta_{i,j}.
\end{align*}
Condition \eqref{item: main theorem, stare su tutti i polinomi} of \Cref{corollary: risultato generale di caratterizzazioni delle soluzioni di un'ODE} says  that the path $X=(t,f(t),f'(t))$ satisfies
\begin{align*}
0&=
\sum_{j=1}^r a_{i,j}\langle \sigma(X),\word{(1+j)}w \rangle- \langle \sigma(X), \word{(1+r+i)}w\rangle
\end{align*}
 for every $i\in\{1,\dots,r\}$ and for every $w\in\TT(\R^d
)$. The converse implication follows from \Cref{corollary: risultato generale di caratterizzazioni delle soluzioni di un'ODE} in the same way.
\end{proof}

 %

 \subsection{Hamiltonian fields}
Let $r=2s$ be even. A special role in classical mechanics is played by vectors fields on the space $\R^r=\R^s\oplus\R^{s\ast}$, that is, global sections 
\begin{equation*}
    \Phi\in\Gamma(\R^s\oplus\R^{s\ast},T(\R^s\oplus\R^{s\ast})).
\end{equation*}
Physically, this corresponds to replacing the space $\R^s\oplus\R^s=T(\R^s)$ of positions and velocities with the space $\R^s\oplus\R^{s\ast}=T^*(\R^s)$ of positions and momenta. Accordingly, we will split the coordinates on $\R^r$ into two subsets: the positions $(x_1,x_2,\ldots,x_s)$ and the so--called conjugate momenta $(p_1,p_2,\ldots,p_s)$. In view of such a splitting, the components of the vector field $\Phi$ (see \cref{rem: integral curve as system of ode}) 
becomes
 \begin{equation*}\label{eqHamVF-proto}
     \Phi=\sum_{i=1}^s\Phi_i\frac{\partial}{\partial x_i}+\Psi_i\frac{\partial}{\partial p_i}\, ,
 \end{equation*}
for suitable functions $\Phi_i(x_1,\dots,x_s,p_1,\dots,p_s)$ and $\Psi_i(x_1,\dots,x_s,p_1,\dots,p_s)$ on $\R^r$.
\begin{definition}
Let $r=2s$ be even. A vector field $\Phi_H\in\Gamma(\R^s\oplus\R^{s\ast},T(\R^s\oplus\R^{s\ast}))$ on $\R^r$ is \textit{Hamiltonian} if there exists a differentiable function $H(x_1,\dots,x_s,p_1,\dots,p_s)$ on $\R^r$ such that 
   \begin{equation*}\label{eqHamVF}
     \Phi_H=\sum_{i=1}^s\frac{\partial H}{\partial p_i}\frac{\partial}{\partial x_i}-\frac{\partial H}{\partial x_i}\frac{\partial}{\partial p_i}.
 \end{equation*}
 \end{definition}
In classical mechanics, usually $H(\x,\p)=V(\x)+\frac{1}{2}Q(\p)\, ,$ where the potential $V(\x)$ is a function on $\R^s$ and the energy $\tfrac{1}{2}Q(\p)$ is a quadratic form on $\R^{s\ast}$. In this case there exists a symmetric $s\times s$ matrix $A$ such that
\begin{equation*}
H(\x,\p)=V(\x)+ \frac{1}{2} \p^\top A\p
\end{equation*}
and the Hamiltonian vector field is
$$
\Phi_H=\sum_{i=1}^s (A\p)_i\frac{\partial}{\partial x_i}-\left( \frac{\partial V}{\partial x_i}(\x)\right)\frac{\partial}{\partial p_i}.
$$
Therefore, for all $i\in\{1,\dots,s\}$, the functions defining the system of ODEs of the integral curves of $\Phi_H$ are
\begin{eqnarray*}
    F_i(t,\x,\p,\y,\q)&=&   A_i\cdot\p -y_i, \\
    F_{s+i}(t,\x,\p,\y,\q)&=& -\frac{\partial V}{\partial x_i}(\x)-q_i,
\end{eqnarray*} 
where $A_i$ is the $i$--th column (or row) of $A$.

\begin{corollary}\label{corollary: Hamiltonian paths}
Let $\vv\in\R^s$ and let $A=(a_{i,j})\in \R^{s\times s}$ be a symmetric matrix. Let $r=2s$ and consider the function $H(\x,\p)=\vv\cdot\x+\tfrac{1}{2}\p^\top\cdot A\cdot\p$ on $\R^r$  fulfilling $\x(0)=\x_0$ and $\p(0)=\p_0$. Then the curve $[a,b]\to\R^r$ given by $t\mapsto(\x(t),\p(t))$ 
is the  integral curve of the Hamiltonian vector field $\Phi_H$ on $\R^r$ fulfilling $(\x(0),\p(0))=(\x_0,\p_0)$ if and only if  the path $X(t)=(t,\x(t),\p(t),\x'(t),\p'(t))$ satisfies 
\begin{align*}
& X(0)=(0,\x_0,\p_0,A\cdot\x_0,-\vv)\\
& \sum_{h=1}^s \langle \sigma(X),a_{i,h}(\word{1}+s+\word{h})w \rangle=-\sum_{h=1}^s\langle \sigma(X),(\word{1}+\word{r}+\word{h})w \rangle\\
& \left\langle\sigma(X), (\word{1}+3s+\word{i}) w\right\rangle=0
\end{align*}
for every $i\in\{1,\dots,s\}$ and for every word $w\in \TT(\R^d)$.
\end{corollary}

\begin{proof}
For every $i\in\{1,\dots,s\}$, define
 \begin{align*}\label{eqDefEffeICampHam2}
 F_i(t,\x,\p,\y,\q)&=  (A\p)_i -y_i\\
    F_{s+i}(t,\x,\p,\y,\q)&= -v_i-q_i
\end{align*} 
and let $f$ be a solution of $$
\begin{cases}
F_1(t,\x(t),\p(t),{\x'}(t),\p'(t))=0\\
\vdots \\
F_r(t,\x(t),\p(t),{\x'}(t),\p'(t))=0.
\end{cases}
$$
In orderto apply  condition \eqref{item: main theorem, stare su tutti i polinomi} of \Cref{corollary: risultato generale di caratterizzazioni delle soluzioni di un'ODE} we have to look at the expression of the $\widetilde{F}_i$'s.
For every $i\in\{1,\dots,s\}$ we have
\begin{align*}
    \widetilde{F}_i(t,\x,\p,\y,\q)&=F_i((t,\x,\p,\y,\q)-X(0))-F_i(X(0))\\
    &=\sum_{l=1}^s a_{i,l}(p_l+f_l(0))-(y_i-f'_i(0))-\left( \sum_{l=1}^s a_{i,l}f_l(0)- f'_i(0)\right)\\
    &=\sum_{l=1}^s a_{i,l}p_l-y_i=F_i(t,\x,\p,\y,\q)
\end{align*}
and
\begin{align*}
 \widetilde{F}_{s+i}(t,\x,\p,\y,\q)&=F_{s+i}((t,\x,\p,\y,\q)-X(0))-F_{s+i}(X(0))\\    
 &=-v_i-(q_i-f'_{s+i}(0))-(-v_i-f'_{s+i}(0))\\
 &=-q_i+2f'_{s+i}(0)\, .
\end{align*}
Moreover
\begin{align*}
    \frac{\partial \widetilde{F}_i}{\partial t}&=0,\quad   \frac{\partial \widetilde{F}_i}{\partial x_j}=0,\quad \frac{\partial \widetilde{F}_i}{\partial p_j}=a_{i,j},\quad \frac{\partial \widetilde{F}_i}{\partial y_j}=-1,\quad \frac{\partial \widetilde{F}_i}{\partial q_j}=0\\
     \frac{\partial \widetilde{F}_{s+i}}{\partial t}&=0,\quad   \frac{\partial \widetilde{F}_{s+i}}{\partial x_j}=0,\quad \frac{\partial \widetilde{F}_{s+i}}{\partial p_j}=0,\quad \frac{\partial \widetilde{F}_{s+i}}{\partial y_j}=0,\quad \frac{\partial \widetilde{F}_{s+i}}{\partial q_j}=-\delta_{i,j}.
\end{align*}
Hence, condition \eqref{item: main theorem, stare su tutti i polinomi} of \Cref{corollary: risultato generale di caratterizzazioni delle soluzioni di un'ODE} reads
\begin{align*}
\sum_{h=1}^s\langle \sigma(X), a_{i,h}(\word{1}+s+\word{h}) w+(\word{1}+\word{r}+\word{h})w  \rangle=0
\hbox{ and }\left\langle\sigma(X), -(\word{1}+3s+\word{i}) w\right\rangle=0
\end{align*}
for all words $w\in \TT(\RR^d)$.
\end{proof}

 \section{Conclusions and perspectives}\label{section: conclusions}
The results we presented hardly cover all the possible applications of signature tensors in this area. We conclude the paper with an outline of what we think could be the next steps in this direction.

\begin{enumerate}
\item Corollaries \ref{coroll: caratterizzazione dei paths che stanno su una varietà} and \ref{theorem: caratterizzazione dei paths che stanno su un'ipersuperficie algebrica} deal with algebraic varieties, that is sets defined by polynomial equations. It is natural to wonder if there are analogous characterizations for manifolds defined by analytic equations. This would allow to generalize \Cref{corollary: risultato generale di caratterizzazioni delle soluzioni di un'ODE} to systems of differential equations that are not necessarily polynomial.
\item Little is known on the rank of signature tensors, with the exception of piecewise linear paths \cite[Section 3]{GalSan}. If $X$ is Legendrian, or if $X$ lies on a given algebraic variety, could we bound $\rk(\sigma^{(k)}(X))$?
\item One of the ideas in \cite{AFS19} is to choose a class $\mathcal{P}$ of paths and the order $k\in\N$ of the signature and study the geometry of $\{\sigma^{(k)}(X)\mid X\in \mathcal{P}\}\subseteq (\R^d)^{\otimes k}$. It would be interesting to do it for the class of holonomic paths, or for the class of paths lying on a given variety.
\item In \Cref{remark: delle nostre equazioni non ne bastano un numero finito} we stress that \Cref{theorem: caratterizzazione dei paths che stanno su un'ipersuperficie algebrica} does not hold if we consider only finitely many equations. This suggests that our method could be useful not just in a qualitative way, but also in a quantitative way. Even if a path $X$ does not lie on a variety $V$, we may define $X$ to be close to be contained in $V$ if $\sigma(X)$ satisfies equations \eqref{eqCond2ThMAIN} up to a certain level of the signature.
\item Finally, since we were interested in solution of Cauchy problems, we restricted our attention to continuous and piecewise differentiable paths, as in \Cref{definition: paths of class Ck}. However, signature tensors are defined for much larger classes of paths, and the uniqueness \Cref{theorem: chen uniqueness} holds even for rough paths. It would be very interesting to generalize our approach in \Cref{section: ODE} to characterize solutions of stochastic differential equations, rather than just ordinary differential equations.
\end{enumerate}

\bibliographystyle{alpha}
\bibliography{References.bib}

\end{document}